\documentclass[12pt]{article}
\usepackage{amsmath,latexsym,amssymb,eucal}
\newtheorem{thm}{Theorem}[section]
\newtheorem{lem}[thm]{Lemma}
\newtheorem{prop}[thm]{Proposition}
\newtheorem{df}[thm]{Definition}
\newtheorem{cor}[thm]{Corollary}
\newcommand{\Ad}{\mathrm{Ad}\,}

\newcommand{\Mor}{\mathrm{Mor}\,}
\newcommand{\Aut}{\mathrm{Aut}\,}

\newcommand{\id}{\mathrm{id}}

\newcommand{\hG}{\hat{G}}
\newcommand{\iG}{\mathrm{Irr}(G)}
\newcommand{\rG}{\mathrm{Rep}(G)}
\newcommand{\cR}{\mathcal{R}}
\newcommand{\tU}{\tilde{U}}
\newcommand{\tV}{\tilde{V}}
 \date{}
 \author{MASUDA Toshihiko\footnote{e-mail address masuda@math.kyushu-u.ac.jp}, \\
 Graduate School of Mathematics, Kyushu University,\\
 6-10-1 Hakozaki, Fukuoka, 812-8581, JAPAN}
 
\begin{document}

 \title{Classification of actions of  duals of finite groups \\ on the AFD
 factor of type II$_1$}
 \maketitle

 \begin{abstract}
 We will show the uniqueness of outer coactions of finite groups on the
 AFD factor of type II$_1$ along the arguments by Connes, Jones and
 Ocneanu. Namely, we construct the infinite tensor product type action,
 adopt it as the model action, and prove that any outer coaction is
  conjugate to the model action.
 \end{abstract}

 \section{Introduction}
 In the theory of operator algebras, the study of automorphisms is one of
 the most important topics. Especially,
 much progress has been made on 
 classification of automorphisms and 
 group actions on injective factors since 
 fundamental works of A. Connes.
 In \cite{Con-auto} and \cite{Co-peri}, A. Connes succeeded in
 classifying automorphisms of the approximately finite dimensional (AFD) 
 factor of type II$_1$ up to outer
 conjugacy. The first generalization of Connnes' results was made by
 V. F. R. Jones in \cite{J-act}, where he classified actions of finite
 groups on the AFD
 factor of type II$_1$. Soon after Jones' theory,
 A. Ocneanu classified actions of discrete amenable groups on the AFD
 factor of type II$_1$. One of extension of their results is
 analysis (or classification) of actions of dual object of groups, 
 i.e., coaction of groups,
 which will be
 useful to understand actions of compact groups. 
 (See \cite{Nk-Tak} on basic of coactions.)

 In this paper, we give the classification
 theorem for outer coactions of finite groups. Here we have to remark
 that this follows indirectly from the works  cited in above. 
 In fact, the uniqueness of outer coactions of finite groups follows from
 \cite{J-act}, since every outer coaction of a finite group is dual to
 some outer (usual) action.
 Nowadays, this also follows from the general theory of Popa's
 classification of subfactors \cite{Po-amen} (also see \cite{PoWa}).
 However in these approach, one does not handle coactions directly.
 Hence in this paper, we present the direct
 approach for classification theorem of outer coactions of finite groups
 on the AFD factor of type II$_1$, which can be generalized to finite
 dimensional Kac algebras. Our argument is similar to
 Connes-Jones-Ocneanu theory. We construct the model
 action on the AFD factor of type II$_1$, prove several cohomology
 vanishing theorem, and compare a given action to the model action. 
 The main technique in this arguments is the ultraproduct and the
 central sequence algebra. Unfortunately, 
 coactions do not necessary induce coactions on the central sequence
 algebra unlike the usual group action case. Hence we have to modify actions to
 apply the ultraproduct technique to handle with coactions, and this is
 one of the important point in our theory.

 Here we have another formulation to treat actions of duals of (finite) groups
 other than coactions due to Roberts in \cite{Rob-act}. His
 approach is essentially equivalent to coactions. However it is
 convenient (at least for the author)
 to regard coactions as the Roberts type actions, 
 which we often call actions of of finite group duals.
 Hence in this paper, we present our main theorem as the uniqueness of Roberts
 type actions of finite group duals.

 This paper is organized as follows. In \S \ref{sec:pre}, we prepare 
 notations used in this paper, and discuss the Roberts type actions.   
 In \S \ref{sec:model}, we construct the infinite tensor product type
 action, which we call the model action. In \S \ref{sec:tech}, we collect
 some technical lemmas, which is  necessary  to treat actions on the
 ultraproduct algebra in \S \ref{sec:ultra}. In \S \ref{sec:coho}, we
 show three kinds of cohomology vanishing theorem, which are important
 tools for analysis of actions. The contents in \S \ref{sec:ultra} and \S \ref{sec:class}
 are central in this paper. We discuss actions
 on the ultra product algebra, and the central sequence algebra. By means
 of cohomology vanishing, we construct the piece of the model action, and
 complete classification. In appendix, we present the Roberts type action
 approach for (twisted) crossed product construction for actions of
 finite group duals.

 \section{Preliminaries and Notations}\label{sec:pre}
 \subsection{Notations on duals of finite groups}\label{sub:inter}
 Throughout this paper, we always assume that $G$ is a finite group. 
 Let $\rG$ and $\iG$  be the collection of all finite dimensional unitary
 representations, and irreducible unitary representations of $G$ respectively.  
 We denote the trivial representation by $\mathbf{1}$.
 We fix representative elements of $\iG/\!\sim$, where $\sim$ means a usual
 unitary equivalence,   
 and denote by $\hG$, and assume $\mathbf{1}\in \hG$.

 Let $d\pi:=\dim H_\pi$ be the  dimension of $\pi\in \rG$.
 For $\pi,\rho\in \rG$, we
 denote the intertwiner space between $\sigma$ and $\pi$ by
 $(\sigma,\pi):=\{T\in B(H_\sigma,H_\pi)\mid T\sigma(g)=\pi(g)T, g\in
 G\}$. If $\sigma$ is irreducible, 
 $(\sigma,\pi)$ becomes a Hilbert space
 with an inner product 
 $\langle T,S \rangle 1:=S^*T$. 

 Let $\pi,\rho\in \hG$, and $\pi\otimes\rho\cong \oplus_{\sigma\in \hG} 
 N_{\pi\rho}^\sigma \sigma$ be the irreducible decomposition, where
 $N_{\pi\rho}^\sigma$ is a multiplicity.
 Fix an orthonormal basis
 $\{T_{\pi,\rho}^{\sigma,e}\}_{e=1}^{N_{\pi\rho}^\sigma}\subset
 (\sigma,\pi\otimes\rho)$. 
 Then we have
 $T_{\pi,\rho}^{\sigma,e*}T_{\pi,\rho}^{\xi,f}=\delta_{\sigma,\xi}\delta_{e,f}1_\sigma$, and 
 $\sum_{\sigma,e}T_{\pi,\rho}^{\sigma,e}T_{\pi,\rho}^{\sigma,e*}=1_{\pi\otimes
 \rho}$, where $1_\sigma\in (\sigma,\sigma)$ is an identity.
 Hence we have $\pi(g)\otimes\rho(g)=
 \sum_{\sigma,e}T_{\pi,\rho}^{\sigma,e}\sigma(g)T_{\pi,\rho}^{\sigma,e*}$ especially.

 Let $\{v(\pi)\}_{\pi\in \hG}$ with $v(\pi)\in
 M_{d\pi}(\mathbf{C})$. Then
 $\sum_{e}T_{\pi,\rho}^{\sigma,e}v(\sigma)T_{\pi,\rho}^{\sigma,e*}$
 does not depend on the choice of $\{T_{\pi,\rho}^{\sigma,e}\}\subset
 (\sigma,\pi\otimes \rho)$.

 In a similar way, one can easily see 

 $$\sum_{\eta,a,b} (T_{\pi,\rho}^{\eta,a}\otimes 1_\sigma)T_{\eta,\sigma}^{\xi,b}v(\xi)
 T_{\eta,\sigma}^{\xi,b*}(T_{\pi,\rho}^{\eta,a*}\otimes 1_\sigma)=
 \sum_{\zeta,c,d}(1_\pi\otimes T_{\rho,\sigma}^{\zeta,c})T_{\pi,\zeta}^{\xi,d}v(\xi)
 T_{\pi,\zeta}^{\xi,d*}(1_\pi\otimes T_{\rho,\sigma}^{\zeta,c*})
 $$
 since both $\{(T_{\pi,\rho}^{\eta,a}\otimes
 1_\sigma)T_{\eta,\sigma}^{\xi,b}\}_{\eta,a,b}$ and 
 $\{(1_\pi\otimes
 T_{\rho,\sigma}^{\zeta,c})T_{\pi,\zeta}^{\xi,d}\}_{\zeta,c,d}$ are
 orthonormal basis for $(\xi,\pi\otimes \rho\otimes\sigma)$. \\
 \noindent
 \textbf{Remark.}
 Assume that $\{v(\pi)\}_{\pi\in \hG}$, $v(\pi)\in A\otimes B(H_\pi)$, is
 given for some vector space $A$. 
 We can extend $v(\pi)$ for a general $\pi\in \rG$ as follows. Let $\pi\cong\oplus_i
 \sigma^i$, $\sigma^i\in\hG$, be an irreducible decomposition, and fix
 $T^i\in(\sigma^i,\pi)$ with $T^{i*}T^j=\delta_{i,j}$ and
 $\sum_{i}T^iT^{i*}=1$. Define
 $v(\pi)=\sum_{i}T^iv(\sigma^i)T^{i*}\in A\otimes B(H_\pi)$. 
 Then $v(\pi)$ is well-defined,
 i.e., it is independent from the choice of $\{T^i\}$,  
 and satisfies $v(\pi)T=Tv(\sigma)$ for $T\in (\sigma,\pi)$.
 In this notation, 
 the contents of the previous paragraph is written as 
 $v((\pi\otimes \rho)\otimes \sigma)=v(\pi\otimes (\rho\otimes \sigma))$,
 for example. 

 Let
 $\{e^\pi_i\}_{i=1}^{d\pi}$ be an orthonormal basis 
 for $H_\pi$, and fix it.
 Let us express
 $T_{\pi,\rho}^{\sigma,e}=(T_{\pi_i,\rho_k}^{\sigma_m,e})$ as a matrix form. 
 Then we can write
 $T_{\pi,\rho}^{\sigma,e*}T_{\pi,\rho}^{\xi,f}=\delta_{\sigma,\xi}\delta_{e,f}1_{\sigma}$
 and
 $\sum_{\sigma,e}T_{\pi,\rho}^{\sigma,e}T_{\pi,\rho}^{\sigma,e*}=1_{\pi\otimes
 \rho}$ by matrix coefficients as 
 $$\sum_{i,k}\overline{T_{\pi_i,\rho_k}^{\sigma_m,e}}T_{\pi_i,\rho_k}^{\xi_n,f}=\delta_{\sigma,\xi}
 \delta_{e,f}\delta_{m,n},$$ and 
 $$\sum_{\sigma,m,e}T_{\pi_i,\rho_k}^{\sigma_m,e}\overline{T_{\pi_j,\rho_l}^{\sigma_m,e}}
 =\delta_{i,j}\delta_{k,l}.$$

 Let $T_{\pi,\bar{\pi}}^\mathbf{1}\in (\mathbf{1},\pi\otimes\bar{\pi})$ 
 be an isometry given by
 $T_{\pi,\bar{\pi}}^{\mathbf{1}}1=\frac{1}{\sqrt{d\pi}}\sum_{i}e^\pi_i\otimes
 e^{\bar{\pi}}_i$, and fix it.
 It is easy to see
 $T_{\pi_i,\bar{\pi}_j}^\mathbf{1}=
 \frac{\delta_{i,j}}{\sqrt{d\pi}}$. 
 Since
 $T_{\pi,\bar{\pi}}^{\mathbf{1}*}T_{\pi,\bar{\pi}}^{\rho,e}=\delta_{\mathbf{1},\rho}$, we
 have $\sum_kT_{\pi_k,\bar{\pi}_k}^{\rho_l,e}=\sqrt{d\pi}\delta_{\mathbf{1},\rho}$

 Set
 $\tilde{T}_{\bar{\pi},\sigma}^{\rho,e}:=
 \frac{\sqrt{d\rho d\pi}}{\sqrt{d\sigma}}(1_{\bar{\pi}}\otimes T_{\pi,\rho}^{\sigma,e*})
 (T_{\bar{\pi},\pi}^\mathbf{1}\otimes 1_\rho)\in (\rho,\bar{\pi}\otimes\sigma)$. 
 Then $\{\tilde{T}_{\bar{\pi},\sigma}^{\rho,e}\}$
 is an
 orthonormal basis for $(\rho,\bar{\pi}\otimes \sigma)$.
 It is easy to see
 $\tilde{T}_{\bar{\pi}_i\sigma_m}^{\rho_k,e}=\sqrt{\frac{d\rho}{d\sigma}}
 \overline{T_{\pi_i,\rho_k}^{\sigma_m,e}}$.
 As a consequence we have 
 $$\sum_{\sigma,m,n,e}
 T_{\pi_i,\rho_k}^{\sigma_m,e}v(\sigma)_{m,n}\overline{T_{\pi_j,\rho_l}^{\sigma_n,e}}
 =\sum_{\sigma,m,n,e}\frac{d\sigma}{d\rho}
 \overline{T_{\bar{\pi}_i,\sigma_m}^{\rho_{k},e}}v(\sigma)_{m,n}
 {T_{\bar{\pi}_j,\sigma_n}^{\rho_l,e}}
 $$
 for example.

 \subsection{Coactions and Roberts type actions}\label{sub:roberts}

 Let $A, B$ be von Neumann algebras.
 We denote the set of unital $*$-homomorphisms from $A$ to $B$ by $\Mor(A,B)$.

 Let $u_g$ be the (right) regular representation of $G$, and $R(G):=\{u_g\}''$
 the group algebra. The coproduct $\Delta$ of $R(G)$ is given by
 $\Delta(u_g)=u_g\otimes u_g$.

 For simplicity, we denote $1_M\otimes T\in M\otimes
 B(H_\pi,H_\rho)$, $T\in B(H_\pi,H_\rho)$,  by $T$.

 \begin{df}\label{df:actiondual}
  $(1)$ Let $M$ be a von Neumann algebra. We say
  $\alpha=\{\alpha_\pi\}_{\pi\in \rG}$ is an action  of $\rG$ if
  $\alpha_\pi\in \Mor(M,M\otimes B(H_\pi))$, and following  hold. \\
 $\mathrm{(1a)}$ $\alpha_\mathbf{1}=\id_M$. \\
  $\mathrm{(1b)}$ $\alpha_\pi(x)T=T\alpha_\sigma(x)$ for any $T\in (\sigma,\pi)$. \\
  $\mathrm{(1c)}$ $\alpha_\pi\otimes \id_\sigma\circ\alpha_\sigma=\alpha_{\pi\otimes \sigma}$. \\
   $(2)$ We say $\alpha=\{\alpha_\pi\}_{\pi\in \iG}$ is an action of $\iG$ if 
 $\alpha_\pi\in \Mor(M,M\otimes B(H_\pi))$ and 
  following holds. \\
 $\mathrm{(2a)}$ $\alpha_\mathbf{1}=\id_M$. \\
  $\mathrm{(2b)}$ $\alpha_\pi\otimes \id_\rho\circ\alpha_\rho(x)T=T\alpha_{\sigma}(x)$ for 
 any $T\in(\sigma,\pi\otimes \rho)$. \\ 
 $(3)$ We say $\alpha=\{\alpha_\pi\}_{\pi\in \hG}$ an action of $\hG$ if 
 $\alpha_\pi\in \Mor(M,M\otimes B(H_\pi))$ and we have  the following. \\
 $\mathrm{(3a)}$ $\alpha_\mathbf{1}=\id_M$. \\
  $\mathrm{(3b)}$ $\alpha_\pi\otimes \id_\rho\circ\alpha_\rho(x)T=T\alpha_{\sigma}(x)$ for 
 any $T\in(\sigma,\pi\otimes \rho)$. \\ 
 \end{df}

 If an action $\alpha$ of $\hG$ is given, then it is a routine work to
 extend $\alpha$ to those of $\iG$ and $\rG$. 
 Hence in this paper, we do not distinguish these notions.
 When $M$ is properly infinite, it is not difficult to see Definition
 \ref{df:actiondual} is reduced to that of Roberts action.

 We remark that $\alpha_\pi$ is automatically injective. Suppose
 $\alpha_\pi(x)=0$. Then we have $0=T_{\bar{\pi},\pi}^{\mathbf{1}*}\alpha_{\bar{\pi}}\otimes
 \id_\pi\circ \alpha_\pi(x)T_{\bar{\pi},\pi}^\mathbf{1}=
 T_{\bar{\pi},\pi}^{\mathbf{1}*}T_{\bar{\pi},\pi}^\mathbf{1}\alpha_\mathbf{1}(x)=x$.

 Let $\{e_{ij}^\pi\}$ be a system of matrix units for $B(H_\pi)$. 
  Then
 $\alpha_\pi(x)$ is decomposed as
 $\alpha_{\pi}(x)=\sum_{i,j}\alpha_\pi(x)_{ij}\otimes e_{ij}^\pi$. The $*$-homomorphism
 property of $\alpha_\pi$ implies
 $\sum_{k}\alpha_\pi(x)_{ik}\alpha_\pi(y)_{kj}=\alpha_\pi(xy)_{ij}$ and
 $(\alpha_\pi(x)_{ij})^*=\alpha(x^*)_{ji} $. 

 The group algebra $R(G)$ can be decomposed as $R(G)=\bigoplus_{\pi\in
 \hG}B(H_\pi)$, 
 and $e^\pi_{ij}=d\pi/|G|\sum_g\pi(g)_{ij}u_g$ gives a matrix unit for
 $B(H_\pi)$. Then $\alpha\in
 \Mor(M,M\otimes R(G))$ can be decomposed as
 $\alpha(x)=\sum_{\pi}\alpha_\pi(x)_{ij}(x)\otimes e^\pi_{ij}$, and we
 get $\alpha_\pi\in\Mor(M, M\otimes B(H_\pi) )$.
 One can verify that $\alpha$ is a coaction, i.e.,
 $\alpha$ is injective and satisfies $(\alpha\otimes
 \id)\circ\alpha=(\id\otimes\Delta)\circ \alpha$. 
  if and only if
 $\{\alpha_\pi\}$ is an action of $\hG$ in the sense of Definition \ref{df:actiondual}.

 \begin{df}\label{df:fix}
  Let $\alpha$ be an action of $\hG$ on $M$. The fixed point algebra
  $M^\alpha$ is defined as $M^\alpha:=\{a\in M\mid \alpha_\pi(a)=a\otimes 1_\pi \,
 \mbox{ for any }\pi\in \hG\}$.
 \end{df}

 If $K\subset M^\alpha $, then we say $\alpha$ is trivial on $K$, and
 often write as $\alpha_\pi=\id$ on $K$. 
 Let $K\subset M$ be a von Neumann subalgebra, on which $\alpha$ acts trivially.
 Then it is easily seen that $\alpha_\pi(K'\cap M)\subset
 (K'\cap M)\otimes B(H_\pi)$, and $\alpha$ is an action on $K'\cap M$.  
 Note that 
 even if we have $\alpha_\pi(K)\subset K\otimes B(H_\pi)$, $\alpha$
 does not induce an action on $K'\cap M$ in general
 unlike the usual group action case.

 Let $\alpha$ be an action of $\hG$ on $M$, and $N$ be another von
 Neumann algebra. Then $\alpha'_\pi(x):=\sum_{i,j}\alpha_\pi(x)_{ij}\otimes
 1_N\otimes e^\pi_{ij}$ is an action of $\hG$ on $M\otimes N$, which we
 denote by $\alpha\otimes \id_N$ for simplicity.

 \subsection{Crossed product construction by Roberts type action}\label{sub:cross}
 Let $\alpha$ be a coaction of $G$ on $M$. 
 The crossed
 product $M\rtimes_\alpha \hG$ is defined as $\alpha(M)\vee \mathbf{C}\otimes
 \ell^\infty(G)\subset M\otimes B(\ell^2(G))$. We discuss the crossed
 product construction from the point of view of the Roberts type action.
 (Also see Appendix.)

 We begin with the following definition.
 \begin{df}\label{df:rep}
  Let $M$ be a von Neumann algebra. We say $\{U_\pi\}_{\pi\in \iG}$ is a
  $($unitary$)$ representation of $\iG$ in $M$ if we have the following. \\
  $(1)$ $U_\pi \in U(M\otimes B(H_\pi))$, $U_\mathbf{1}=1$. \\
  $(2)$ Let $F_{\pi,\sigma}\in B(H_\pi\otimes H_\rho,H_\rho\otimes H_\pi)$ 
 be a flip map. Set $U_\pi^{12}:=U_\pi\otimes 1_\rho$, and
  $U_\rho^{13}:=F_{\rho,\pi}(U_\rho\otimes 1_\pi)F_{\pi,\rho}$. 
 Then $U_\pi^{12}U_\rho^{13} T=T U_\sigma$ for any $T\in (\sigma,\pi\otimes \rho)$.
 \end{df}

 If we represent $U_\pi$ and $T$ as $U_\pi=(U_{\pi_{ij}})_{1\leq i,j \leq d\pi}$ and
 $T=(T_{i,k}^m)_{1\leq i\leq d\pi, 1\leq k \leq d\rho}^{1\leq m\leq
 d\sigma}$ 
 respectively by matrix elements,
 then Definition \ref{df:rep}(2) is written as
 $\sum_{j,l}U_{\pi_{ij}}U_{\rho_{kl}}T_{j,k}^n
 =\sum_mT_{i,k}^mU_{\sigma_{mn}}$. 

 \begin{lem}\label{lem:reprel}
 Let $\{U_\pi\}$ be a unitary representation of $\hG$. Then we have
  $U_{\pi_{ij}}^*=U_{\bar{\pi}_{ij}}$, and $[U_{\pi_{ij}},U_{\rho_{kl}}]=0$.
 \end{lem}
 \textbf{Proof.} Since we have
 $\sum_{j,l}U_{\pi_{ij}}U_{\bar{\pi}_{kl}}T_{\pi_j,\bar{\pi}_l}^\mathbf{1}=
 T_{\pi_i,\bar{\pi}_k}^\mathbf{1}U_{\mathbf{1}}$,
 $\sum_{j}U_{\pi_{ij}}U_{\bar{\pi}_{kj}}=\delta_{ik}$ holds. This implies
 $U_\pi{}^tU_{\bar{\pi}}=1$, and hence
 $U_{\pi}^*={}^tU_{\bar{\pi}}$. Thus we get $U_{\pi_{ij}}^*=U_{\bar{\pi}_{ij}}$. 

 We will verify the second statement. 
 Since $U_\pi$ is a representation, we have
 $U_\pi^{12}U_\rho^{13}=\sum_{\sigma,e}T_{\pi,\rho}^{\sigma,e }U_\sigma 
 T_{\pi,\rho}^{\sigma,e* }$.
 Then we get
 $$F_{\pi,\rho}(U_\pi^{12} U_\rho^{13}) F_{\rho,\pi}=
 \sum_{i,j,k,l}U_{\pi_{ij}}U_{\rho_{kl}}\otimes
 e_{kl}^\rho\otimes e_{ij}^\pi
 =\sum_{\sigma,e}(F_{\pi,\rho}T_{\pi,\rho}^{\sigma,e })U_\sigma 
 (F_{\pi,\rho}T_{\pi,\rho}^{\sigma,e})^*.$$ 

 On the other hand, 
 $U_\rho^{12}U_\pi^{13}=\sum_{\sigma,e}F_{\pi,\rho}T_{\pi,\rho}^{\sigma,e }U_\sigma 
 (F_{\pi,\rho}T_{\pi,\rho}^{\sigma,e })^*$ holds, since
 $\{F_{\pi,\rho}T_{\pi,\rho}^{\sigma,e}\}
 \subset
 (\sigma,\rho\otimes \pi)$ is an orthonormal basis. (Note that we use
 $\pi\otimes \rho\sim \rho\otimes \pi$ here.)
 By comparing these, we get  
 $[U_{\pi_{ij}},U_{\rho_{kl}}]=0$. \hfill$\Box$

 Lemma \ref{lem:reprel} shows that $\{U_{\pi_{ij}}\}$ behave like matrix
 coefficients $\{\pi(g)_{ij}\}$.
 Let $U_\pi$ be a representation of $\hG$,
 then it follows immediately that so is $U_\pi^*$, since $[U_{\pi_{ij}},
 U_{\rho_{kl}}]=0$. \\ 

 \noindent
 \textbf{Remark.}
 One can see that $\{U_{\bar{\pi}}^*\}$ is a conjugate representation of $\hG$,
 i.e., $(U_{\bar \pi}^*)^{12}
 (U_{\bar{\rho}}^*)^{13}\overline{T}=\overline{T}U_{\bar{\sigma}}^*$ for
 $T\in (\sigma,\pi\otimes \rho)$,  without using the commutativity of
 $\hG$. \\

 Let $\pi(g)_{ij}$ be a matrix coefficient for $\pi\in \rG$.
 We regard  $\pi(g)_{ij}$ as an element $\pi_{ij} $ in $\ell^\infty(G)$ and set
 $\lambda_{\pi_{ij}}:=1_M\otimes \pi_{ij}$. Then
 $\lambda_\pi=\sum_{i,j}\lambda_{\pi_{ij}}\otimes e^{\pi}_{ij}$ is 
 the unitary representation of $\hG$ in the sense of Definition \ref{df:rep}.

 Since $\ell^\infty(G)=\bigvee\{\pi_{ij}\}$, we have
 $M\rtimes_\alpha \hG=\alpha(M)\vee \{\lambda_{\pi_{ij}}\}$.
 The relation of generators are
 $\sum_{k}\lambda_{\pi_{ik}}x\lambda_{\pi_{jk}}^*
 =\alpha_{\pi}(x)_{ij}$, or equivalently  
 $\lambda_\pi(x\otimes
 1_\pi)\lambda_\pi^*=\alpha_\pi(x)$.
 Here
 we identify $\alpha(x)$ and $x$ as in the usual way.
 A unitary $\lambda_\pi$ plays a roll of
 the implementing unitary in the usual crossed product
 construction. Hence we also call $\lambda_\pi$ the implementing unitary in
 $M\rtimes_\alpha \hG$. 
 We can expand $a \in M\rtimes_\alpha \hG$ as
 $\sum_{\pi,i,j}a_{\pi,i,j}\lambda_{\pi_{ij}}$, $a_{\pi,i,j}\in M$,  
 uniquely.

 \begin{df}\label{df:outer}
  Let $\alpha$ be an action of $\hG$ on $M$. We say $\alpha$ is free if
  there exists no non-zero $a\in M\otimes B(H_\pi)$,  $\mathbf{1}\ne \pi\in \hG$,  
  so that $\alpha_\pi(x)a=a(x\otimes 
  1_\pi)$ for every $x\in M$.
 When $\alpha$ is an action of a factor $M$, then we also say $\alpha$ is
 outer if $\alpha$ is free. 
 \end{df}

 In usual, freeness of a coaction $\alpha$ on a factor $M$ is defined by
 the relative commutant condition
 $M'\cap M\rtimes_\alpha \hG=Z(M)$.
 We see that the usual definition and ours coincide in the
 following proposition.

 \begin{prop}\label{prop:outer}
  Let $\alpha$ be an action of $\hG$ on $M$. Then
  $\alpha$ is  free if and only if  
 $M'\cap M\rtimes_\alpha \hG =Z(M)$. 
 Especially, $M\rtimes_\alpha \hG$ is a factor when $\alpha$ is free, and
  $M$ is a factor.
 \end{prop}
 \noindent
 \textbf{Proof.} Let $a=\sum_{\pi,i,j}a_{\pi,i,j}\lambda_{\pi_{ij}}\in
 M\rtimes_\alpha \hG$. Set $a_\pi:=\sum_{i,j}a_{\pi,j,i}\otimes e^\pi_{ij}\in
 M\otimes B(H_\pi)$. Then it is easy to see $a\in M'\cap M \rtimes_\alpha
 \hG$ if and only if $(x\otimes 1_\pi)a_\pi=a_\pi\alpha_\pi(x)$ for any
 $x\in M$, $\pi \in \hG$. Then it is easily shown that $\alpha$ is free
 if and only if $M'\cap M\rtimes_\alpha \hG = Z(M)$.
 \hfill$\Box$

 In the end of this subsection, we explain the dual action of $G$ on the
 crossed product. Let $\alpha$ be an action of $\hG$ on $M$. Then the dual
 action $\hat{\alpha}$ of $G$ on $M\rtimes_\alpha \hG$ is given by 
 $\hat{\alpha}_g(a)=a$ for $a\in M$, and $\hat{\alpha}_g\otimes
 \id_\pi(\lambda_\pi)=\lambda_\pi \pi(g)$, or equivalently
 $\hat{\alpha}_g(\lambda_{\pi_{ij}})= \sum_k \lambda_{\pi_{ik}}\pi(g)_{kj}$.
 Then 
  it is shown that $\hat{\alpha}$ an action of $G$, and the fixed point
 algebra is  
 $(M\rtimes_\alpha \hG)^{\hat{\alpha}}=M$.

 \subsection{Quantum double construction for  finite group duals}
 In this subsection, we collect definitions and basic properties for
 quantum double construction (also known as the symmetric enveloping
 algebra \cite{Po-symm}, or the 
 Longo-Rehren construction \cite{LR})
  arising from actions of group duals. 
 We will use them in \S \ref{sec:ultra}.  
  We refer \cite[Chapter 12.8, 15.5]{EK-book}, or
 \cite[Appendix A]{M-ext} for details of this topic. 

 For $\pi,\rho\in \rG$, 
 let $\pi\hat{\otimes }\rho$ a representation of $G\times G$ given by
 $\pi\hat{\otimes }\rho(g,h):=\pi(g)\otimes \rho(h)$.
 Let $\alpha$ be an action of $\hG\times \hG$ on $M$. 
 Set $P:=M\rtimes_\alpha (\hG\times \hG)$.  
 Let $\lambda_{\pi
 \hat{\otimes }\rho}$ be an implementing unitary for $\alpha$.  

 \begin{lem}\label{lem:LRunitary}
 Set $w_{\pi_{ij}}:=\sum_k\lambda_{\pi_{ik}\hat{\otimes}\bar{\pi}_{jk}}$. Then
  $w_\pi=(w_{\pi_{ij}})$ is a unitary representation of $\hG$. 
 \end{lem}
 \textbf{Proof.}
 Set $v_{\pi_{ij}}:=\lambda_{\pi_{ij}\hat{\otimes}\mathbf{1}}$, 
 $u_{\pi_{ij}}:=\lambda_{\mathbf{1}\hat{\otimes}\bar{\pi}_{ji}}$.
 Obviously we have
 $w_{\pi_{ij}}=\sum_{k}v_{\pi_{ik}}u_{\pi_{kj}}$ and $[v_{\pi_{ij}}, u_{\rho_{kl}}]=0$. 
 Since $\{\overline{T_{\pi,\rho}^{\sigma,e}}\}\subset (\bar{\sigma},
 \bar{\pi}\otimes \bar{\rho})$ is an orthonormal basis, 
 $u_{\pi}=(u_{\pi_{ij}})_{ij}$ becomes a unitary representation of
 $\hG$ (also see Remark after Lemma \ref{lem:reprel}). 
 Hence 
 \begin{eqnarray*}
  \sum_{j,l}w_{\pi_{ij}}w_{\rho_{kl}}T_{\pi_j,\rho_l}^{\sigma_m,e} &=&
  \sum_{j,l,n,a}v_{\pi_{in}}v_{{\rho_{ka}}}u_{\pi_{nj}}u_{\rho_{aj}}
 T_{\pi_j,\rho_l}^{\sigma_m,e} \\
 & = & \sum_{j,l,n,a,\atop \xi,b,c,p\eta,d,f,q}
 T_{\pi_i,\rho_k}^{\xi_b,p}
 v_{\xi_{bc}}
 \overline{T_{\pi_n,\rho_a}^{\xi_c,p}}
 T_{{\pi}_n,{\rho}_a}^{{\eta}_d,q}
 u_{\eta_{df}}
 \overline{T_{\pi_j,\rho_l}^{\eta_f,q}}
 T_{\pi_j,\rho_l}^{\sigma_m,e} \\ 
 & = & \sum_{\xi,b,c,p\eta,d,f,q}
 T_{\pi_i,\rho_k}^{\xi_b,p}
 \left(\sum_{n,a}
 \overline{T_{\pi_n,\rho_a}^{\xi_c,p}}
 {T_{{\pi}_n,{\rho}_a}^{{\eta}_d,q}}\right)
 v_{\xi_{bc}}u_{\eta_{df}}
 \left(\sum_{j,l}\overline{T_{{\pi}_j,{\rho}_l}^{{\eta}_f,q}}
 T_{\pi_j,\rho_l}^{\sigma_m,e} \right)\\ 
 &=&
 \sum_{b}T_{\pi_i,\rho_k}^{\sigma_b,e}
 \left(\sum_cv_{\sigma_{bc}}u_{\sigma_{cm}}\right) \\
 &=&
 \sum_{b}T_{\pi_i,\rho_k}^{\sigma_b,e}w_{\sigma_{bm}}
 \end{eqnarray*}
 holds. \hfill$\Box$

 \begin{df}
  Set $N:=M\vee \{w_{\pi_{ij}}\}$. We call $M\subset N$ is the quantum double
 for $\alpha$.
 \end{df}

 \noindent
 \textbf{Remark.} In the above definition, we consider an action of
 $\hG\times \hG$ on $M$ directly. However, usual quantum double
 construction is given as follows.
 Let $M$ be a von Neumann algebra, and $\alpha$ be an
 action of $\hG$ on $M$. By the commutativity of $\hG$, 
 $(\alpha_{\bar{\pi}})^\mathrm{opp}$ becomes an action of $\hG$ on
 $M^\mathrm{opp}$. Hence we have an action of $\hG\times \hG$ on
 $M\otimes M^\mathrm{opp}$.
 The rest of construction is same as above. \\

 We embed $G$ into $G\times G$
 by $g\rightarrow (g,g)$. Let $\beta:=\hat{\alpha}$ be the dual action of
 $G\times G$ on $P$. 
 Then it is shown that $N=(M\rtimes_\alpha (\hG\times \hG))^G$.

 For example, we have
 \begin{eqnarray*}
  \beta_{g,g}(w_{\pi_{ij}})&=&\sum_{k}\beta_{g,g}(\lambda_{\pi_{ik}\hat{\otimes}\bar{\pi}_{jk}}) \\
 &=& \sum_{k,l,m}
 \lambda_{\pi_{il}\hat{\otimes}\bar{\pi}_{jm}}\pi(g)_{lk}\overline{\pi(g)_{mk}} \\
 &=&
 \sum_{l,m}
 \lambda_{\pi_{il}\hat{\otimes}\bar{\pi}_{jl}} \\
 &=&w_{\pi_{ij}}.
 \end{eqnarray*}

 If we expand $a\in P$ as 
 $a=\sum a_{\pi_{ij},
 \rho_{kl}}\lambda_{\pi_{ij}\hat{\otimes}\rho_{kl}}$, then 
 $N=(M\rtimes_\alpha (\hG\times \hG))^G$ is verified in a similar way as
 above. We leave the proof to the reader. We remark that $a\in N$ can be
 expand uniquely as  $a=\sum_{\pi,i,j}a_{\pi,i,j}w_{\pi_{ij}}$,
 $a_{\pi,i,j}\in M$, and there exists the canonical conditional
 expectation $E:N\rightarrow M$ given by $E(a)=a_{\mathbf{1}}$.

 \subsection{Main result}

 \begin{df}
  Let $\alpha$ be an action of $\hG$. We say $\{w_\pi\}_{\pi\in \hG}$  a
  $($unitary$)$ 1-cocycle for $\alpha$ if $w_\pi\in U(M\otimes B(H_\pi))$, 
 normalized as $w_\mathbf{1}=1$, and 
 following holds.
 $$(w_\pi\otimes 1_\rho)
 \alpha_\pi\otimes \id_\rho(w_\rho)T=Tw_\sigma,\quad T\in (\sigma, \pi\otimes\rho).$$ 
 A 1-cocycle $\{w_\pi\}$ for $\alpha$ is called a coboundary if there exists a unitary
  $v\in U(M)$ such that $w_\pi=(v^*\otimes 1_\pi)\alpha_\pi(v)$.
 \end{df}
 If we extend $v_\xi$ for $\xi\in \rG$ as in the
 remark in
 \S\ref{sub:inter}, then we have
 $(v_\xi\otimes 1_\eta)\alpha_\xi(v_\eta)=v_{\xi\otimes \eta}$. 
 It is easy to see that $\Ad w_\pi \alpha_\pi$ is an action of $\hG$ for a
 1-cocycle $w_\pi$.

 \begin{df}\label{df:conj}
 Let $\alpha$ and $\beta$ be actions of $\hG$ on $M$. \\
  $(1)$ We say $\alpha$ and
  $\beta$ are conjugate if there exists $\theta\in \Aut(M)$ with 
 $\theta\otimes \id_\pi\circ \alpha_\pi\circ \theta^{-1}=\beta_\pi$ for
  every $\pi\in \hG$. \\
 $(2)$ We say $\alpha$ and $\beta$ are cocycle conjugate if there exists
  a 1-cocycle $\{w_\pi\}$ for $\alpha$, and $\Ad w_\pi \alpha_\pi$ and
  $\beta_\pi$ are conjugate. 
 \end{df}

 Our main purpose is to show the following theorem by the traditional
 Connes-Jones-Ocneanu type approach.

 \begin{thm}\label{thm:main00}
 Let $\cR$ be the AFD factor of type II$_1$.
 Let $\alpha$ and $\beta$ be outer actions of $\hG$ on $\cR$. 
 Then $\alpha$ and $\beta$ are conjugate.
 \end{thm}

 \section{Model action}\label{sec:model}
 In this section,
 we construct an infinite tensor product type action of $\hG$ on $\cR$, 
 which we adopt as the model action.

 It is easy to see the following lemma.
 \begin{lem}\label{lem:reptensor}
  Let $M$, $N$ be von Neumann algebras, and $U_\pi$, $V_\pi$ unitary
  representation of $\hG$ in $M$ and $N$ respectively. We regard $U_\pi$
  and $V_\pi$ as representations of $\hG$ in $M\otimes N$ in the canonical
  way. Then $U_\pi V_\pi$ is also a representation of $\hG$.
 \end{lem}

 To construct the model action, we first construct (the canonical)
 unitary representation of $\iG$ on $M_{|G|}(\mathbf{C})$. Although we
 already discussed it in \S\ref{sub:cross}, we give a slightly different
 approach, which will be useful for our argument.

 Let $\phi$ be the Haar functional for $R(G)$, i.e., $\phi(u_g)=|G|\delta_{e,g}$.
 For $v\in R(G)$, we denote by $v=\oplus v(\pi)$, $v(\pi)\in B(H_\pi)$, via the decomposition
 $R(G)\cong\bigoplus_{\pi\in\hG}B(H_\pi)$. Then we have $\phi(v)=\sum_\pi d\pi
 \mathrm{Tr}_\pi(v(\pi))$, where
 $\mathrm{Tr}_\pi$ be the canonical (non-normalized) trace on $B(H_\pi)$.
 We regard $R(G)$ as a Hilbert space equipped with an inner product
 arising from $\phi$, and denote by $\ell^2(\hG)$. 
 Namely, an inner product on $\ell^2(\hG)$ is given by $\langle
 v,w\rangle=
 \sum_{\pi\in \hG}d\pi \langle v(\pi), w(\pi)\rangle_\pi $ for $v=\oplus 
 v(\pi)$, $w=\oplus w(\pi)$. 
 Here $\langle v(\pi), w(\pi)\rangle_\pi=\mathrm{Tr}_\pi(w(\pi)^*v(\pi))$
 It is easy to see
 $\{d\pi^{-1/2}e^\pi_{ij}\}\subset \ell^2(\hG)$ forms an orthonormal basis
 with respect to this inner product.   

 Set $T_{\rho,\pi_i}^{\sigma,e}\in B(H_\sigma,H_\rho)$ by
 $(T_{\rho,\pi_i}^{\sigma,e})_{\rho_j,\sigma_k}=T_{\rho_j,\pi_i}^{\sigma_k,e}$. 
 \begin{lem}\label{lem:model1}
  Define $\lambda_{\pi_{ij}}\in B(\ell^2(\hG))= M_{|G|}(\mathbf{C})$ by 
 $$(\lambda_{\pi_{ij}}v)(\rho):=
 \sum_{\sigma,e}T_{\bar{\rho},\pi_i}^{\bar{\sigma},e}
 v(\sigma)T_{\bar{\rho},\pi_j}^{\bar{\sigma},e*}
 $$
 and $\lambda_\pi:=\sum_{i,j}\lambda_{\pi_{ij}}\otimes e_{ij}^\pi\in M_{|G|}
 (\mathbf{C})\otimes B(H_\pi)$. Then $\{\lambda_\pi\}$ is a unitary
  representation of $\iG$ on $M_{|G|}(\mathbf{C})$.
 \end{lem}
 \textbf{Proof.}
 We freely use notations and results
 in \S\ref{sub:inter}.
 We first show
 $\lambda_\pi^{12}\lambda_\rho^{13}T_{\pi,\rho}^{\sigma,a}=T_{\pi,\rho}^{\sigma,a}
 \lambda_\sigma$,
 $T_{\pi,\rho}^{\sigma,a}=(T_{\pi_i,\rho_k}^{\sigma_m,a})\in (\sigma,\pi\otimes
 \rho)$, equivalently   
 $\lambda_{\pi_{ij}}\lambda_{\rho_{kl}}T_{\pi_j,\rho_l}^{\sigma_n,a}=\sum_m
 T_{\pi_i,\rho_k}^{\sigma_m,a}\lambda_{\sigma_{mn}}$.
 \begin{eqnarray*}
   \left(\sum_{j,l}\lambda_{\pi_{ij}}\lambda_{\rho_{kl}}
   T_{\pi_j,\rho_l}^{\sigma_n,a}v\right)(\xi) &=& 
   \sum_{j,l,\eta,e}T_{\bar{\xi},\pi_i}^{\bar{\eta},e}
   (\lambda_{\rho_{k,l}} 
   T_{\pi_j,\rho_l}^{\sigma_n,a}v)(\eta)T_{\bar{\xi},\pi_j}^{\bar{\eta},e*}  \\
   &=&
    \sum_{j,l,\eta,e,\zeta,f}T_{\bar{\xi},\pi_i}^{\bar{\eta},e}
	 T_{\bar{\eta},\rho_k}^{\bar{\zeta},f} 
	T_{\pi_j,\rho_l}^{\sigma_n,a}v(\zeta)
 T_{\bar{\eta},\rho_l}^{\bar{\zeta},f*}T_{\bar{\xi},\pi_j}^{\bar{\eta},e*}  \\
 &=& \sum_{j,l,\eta,e,\zeta,f}T_{\pi_i,\rho_k}^{\eta,e}
	 T_{\bar{\xi},\eta}^{\bar{\zeta},f} 
	T_{\pi_j,\rho_l}^{\sigma_n,a}v(\zeta)
 T_{\bar{\xi},\eta}^{\bar{\zeta},f*}T_{\pi_j,\rho_l}^{\bar{\eta},e*}  \\
   &=& \sum_{\eta,\zeta,e,f}
 T_{\pi_i,\rho_k}^{\eta,e} T_{\bar{\xi},\eta}^{\bar{\zeta},f}
 v(\zeta) T_{\bar{\xi},\eta}^{\bar{\zeta},f*}
 (\sum_{j,l}T_{\pi_j,\rho_l}^{\sigma_n,a}
 T_{\pi_j,\rho_l}^{\eta*,e})   \\
   &=& \sum_{\zeta,f} T_{\pi_i,\rho_k}^{\sigma,a}  T_{\bar{\xi},\sigma}^{\bar{\zeta},f}
 v(\zeta) T_{\bar{\xi},\sigma_n}^{\bar{\zeta},f*} \\
   &=& \sum_{m,\zeta,f} T_{\pi_i,\rho_k}^{\sigma_m,a}  T_{\bar{\xi},\sigma_m}^{\bar{\zeta},f}
 v(\zeta) T_{\bar{\xi},\sigma_n}^{\bar{\zeta,f}*} \\
   &=& \sum_{m}(T_{\pi_i,\rho_k}^{\sigma_m,a}  
 \lambda_{\sigma_{m,n}}v)(\xi). 
 \end{eqnarray*}

 It is easy to see that 
 $\sum_k\langle \lambda_{\pi_{ki}}v, \lambda_{\pi_{kj}}w \rangle =\delta_{i,j}$.
 Hence we have
 $\sum_{k}\lambda_{\pi_{ki}}^*\lambda_{\pi_{kj}}=\delta_{i,j}$ and
 consequently $\lambda_\pi^*\lambda_\pi=1$. Thus it suffices to show
 $\lambda_\pi\lambda_\pi^*=1$. 

 Here  we have
 \begin{eqnarray*}
 \sum_{k}\lambda_{\pi_{ik}}\lambda_{\bar{\pi}_{jk}}&=&
  \sum_{k,\rho,l,m,e}T_{\pi_i,\bar{\pi}_j}^{\rho_{l},e}\lambda_{\rho_{lm}}
 \overline{T_{\pi_k,\bar{\pi}_k}^{\rho_{m},e} }\\
 &=&\sqrt{d\pi}
  T_{\pi_i,\bar{\pi}_j}^{\mathbf{1}}\lambda_{\mathbf{1}} \\
 &=& \delta_{i,j}.
 \end{eqnarray*}

 Hence we have $\lambda_\pi{}^t\lambda_{\bar{\pi}}=1$, and
 ${}^t\lambda_{\bar{\pi}}=\lambda_\pi^*$. It follows that
 $\lambda_{\pi_{ij}}^*=\lambda_{\bar{\pi}_{ij}}$ and $\lambda_\pi\lambda_\pi^*=1$.
 \hfill$\Box$

 Let $E=\{e_{\pi_{ij},\rho_{kl}}\}$ be a system of matrix units for
 $B(\ell^2(\hG))\cong 
 M_{|G|}(\mathbf{C})$, that is, $e_{\pi_{ij},\rho_{kl}}$ is a partial
 isometry which sends ${d\rho^{-1/2}}e^\rho_{kl}\in \ell^2(\hG)$ to $d\pi^{-1/2}e^\pi_{ij}$. 
 It is not difficult to see 
 $$ \lambda_{\pi_{ij}}=
 \sum_{\rho,{k, l},\sigma, {m,n},e}\sqrt{\frac{d\rho}{d\sigma}}
 T_{\pi_i,\rho_k}^{\sigma_m,e}\overline{T_{\pi_j,\rho_l}^{\sigma_n,e}}e_{\sigma_{mn},\rho_{kl}}.
 $$
 It follows that
 $\sqrt{d\pi
 d\rho}\lambda_{\pi_{ij}}e_{\mathbf{1},\mathbf{1}}
 \lambda_{\bar{\rho}_{kl}}=e_{\pi_{ij},\rho_{kl}}$ from the above
 expression of $\lambda_\pi$.

 Let $M$ be a von Neumann algebra, and 
 $E=\{e_{\pi_{ij}, \rho_{kl}}\}\subset M$
 a system of matrix units for $B(\ell^2(\hG))$.
 Then we can construct a unitary representation $\lambda_\pi$ of $\hG$ in
 $E''$ by the above formula. 
 In this case, we call $\{\lambda_\pi\}$ a representation of $\hG$
 associated with $E=\{e_{\pi_{ij},\rho_{kl}}\}$.
 When we have to specify $E$, we denote the unitary representation of $\hG$
  associated with $E$ by $\lambda_\pi^E$.

  We define the product type action of $\hG$ on $\cR$. Express
 $\cR=\bigotimes_{n=1}^\infty K_n$, where $K_n$ is a copy of $M_{|G|}(\mathbf{C})$. 
 Let $\lambda^n_{\pi}:=\lambda_\pi^{K_n}$ be a unitary representation of
 $\hG$ on $K_n$, and regard as one on $\cR$.
 Define $\tilde{\lambda}_\pi^1:=\lambda_\pi^1$, and
 $\tilde{\lambda}_\pi^n=
 \tilde{\lambda}_\pi^{n-1}\lambda_\pi^n$.
 Then $\tilde{\lambda}_\pi^n$ is a representation of $\hG$ on $K_1\otimes \cdots
 \otimes K_n$ by Lemma \ref{lem:reptensor}.
 Set $m_\pi^n(x):=\Ad \tilde{\lambda}_\pi^n(x\otimes 1_\pi)$. 
 Since $\tilde{\lambda}_\pi^n$ is a unitary representation of $\hG$,
 $m_\pi^n$ is indeed an action of $\hG$ on $\cR$.
 If $x\in
 \bigotimes_{k=1}^{n-1}K_k$, then 
 $$\Ad \tilde{\lambda}^n_\pi(x\otimes 1_\pi)=\Ad
 \tilde{\lambda}_\pi^{n-1}
 \lambda_\pi^n(x\otimes 1_\pi)=\Ad
 \tilde{\lambda}_\pi^{n-1}(x\otimes 1_\pi)$$
 holds. Hence $\lim\limits_{n\rightarrow \infty}m^n_\pi(x)$ exists for $x\in
 \bigcup_{n=1}^\infty\bigotimes_{k=1}^n K_k$, and so does
 $m_\pi(x)=\lim\limits_{n\rightarrow \infty} m_\pi^n(x)$
 for every $x\in \cR$.

 \begin{df}\label{df:model}
  We call $m=\{m_\pi\}$ the model action for $\hG$.
 \end{df}

 \begin{thm}\label{thm:modelout}
  The model action $m$ is outer.
 \end{thm}
 \textbf{Proof.} Fix $\mathbf{1}\ne \pi \in \hG$. Assume there exists non-zero
 $a\in \cR\otimes B(H_\pi)$ such that 
 $m_\pi(x)a=a(x\otimes 1)$ holds for $x\in \cR$. If $x\in
 \bigotimes_{k=1}^n K_n$, then $(x\otimes
 1)\tilde{\lambda}_\pi^{n*}a=\tilde{\lambda}_\pi^{n*}a(x\otimes 1)$ holds. Hence $a$ is
 expressed as
 $a=\tilde{\lambda}_\pi^nb_{n+1}$, $b_{n+1}=\sum_{ij}b^{n+1}_{ij}\otimes e^\pi_{ij}
 \in \bigotimes_{k=n+1}^\infty K_k\otimes B(H_\pi)$.
 Since we assume $a\ne 0$, there exists $c\in \cR\otimes B(H_\pi)$ with
 $\tau\otimes \mathrm{Tr}_\pi(ca)\ne 0$. We may assume $c$ is of the form $c_1\otimes
 e^\pi_{ij}$, $c_1\in \bigotimes_{k=1}^m K_k $ for some $m$. 
 Then $$\tau\otimes
 \mathrm{Tr}_\pi(ca)=\tau(c_1\lambda_{\pi_{ij}}^{m+1}b_{ji}^{m+2})=
 \tau(c_1\tilde{\lambda}_{\pi_{ij}}^{m+1})\tau(b_{ji}^{m+2})=
 \sum_{l}\tau(c_1\tilde{\lambda}_{\pi_{il}}^m)\tau(\lambda_{\pi_{lj}}^{m+1})\tau(b_{ji}^{m+2})=0
 $$
 holds, and this is a contradiction.
 Hence $a$ must be $0$, and $m$ is an outer action. \hfill$\Box$ 

 \begin{df}\label{df:equmu}
  Let $E=\{e_{\pi_{ij},\rho_{kl}}\}\subset M$ be a system of matrix
  units, and $\lambda_{\pi}^E$ a representation of $\hG$ associated with $E$. 
 Let $\alpha$ be an action of $\hG$ on $M$.
 We say 
 $\{e_{\pi_{ij},\rho_{kl}}\}$ is an $\alpha$-equivariant system of  matrix units if 
 $\alpha_\pi(x)=\Ad \lambda_\pi^E(x\otimes 1)$ for $x\in E$. 
 \end{df}

 The following lemma is easily verified. We leave the proof to the reader.

 \begin{lem}\label{lem:repcocycle}
  Let $\alpha$ be an action of $\hG$ on $M$. \\
 $(1)$ Let $E=\{e_{\pi_{ij},\rho_{kl}}\}$ be an $\alpha$-equivariant
  system of matrix units. Then $\lambda_\pi^{E*}$ is a 1-cocycle for
  $\alpha$, and $\Ad \lambda_\pi^{E*}\alpha_\pi=\id $ on $E$. Hence 
 $\Ad \lambda_\pi^{E*}\alpha_\pi $ induces an action on $E'\cap M$. \\
 $(2)$ Let $M\supset K\cong M_{|G|}(\mathbf{C})$, and suppose $\alpha$ is trivial
  on $K$. Then $\lambda_\pi^K$ is a 1-cocycle for $\alpha$. It follows
  that  $\Ad \lambda_\pi^K\alpha_\pi$ is an action on $M$. 
 \end{lem}

 \section{Technical results}\label{sec:tech}
 In this section, 
 we collect some technical lemmas, whose proof can be found in 
 \cite{Con-auto}, \cite{J-act}, \cite{Ocn-act}.
  In the following, $M$
 is a factor of type II$_1$, and $\tau$ is the unique normalized trace on
 $M$.

 \begin{lem}[{\cite[Lemma 3.2.7]{J-act}}]\label{lem:nearproj}
  Let $f\in M$ be such that $\|f\|\leq 1$, $\|f^2-f\|_2<\delta$ and
  $\|f^*-f\|_2<\delta\leq 1/4$. Then there exists a projection $p\in M$ such that 
 $\|f-p\|_2<6\sqrt[4]{\delta}$ and $\tau(p)=\tau(f)$.
 \end{lem}

 \begin{lem}[{\cite[Lemma 3.2.1]{J-act}}]\label{lem:nearunit}
  Let $u\in M$ be such that $\|u^*u-1\|_2<\delta$. Then there exists a
  unitary $v\in M$ with $\|u-v\|_2<(3+\|u\|)\delta$.
 \end{lem}

 \begin{lem}[{\cite[Proposition 1.1.3]{Con-auto},\cite[Proposition
 7.1]{Ocn-act}}]\label{lem:ultra}
 Let us fix a free ultrafilter $\omega$ over $\mathbf{N}$. \\
 $(1)$ Let $A\in M^\omega$ be a unitary $($resp. projection$)$. Then there
  exists a representing sequence $A=(a_n)$ consisting of unitaries
  $($resp. projections$)$. \\
 $(2)$ Let $V\in M^\omega$ be a partial isometry with $V^*V=E$ and $VV^*=F$.
 Let  $E=(e_n),F=(f_n)\in M^\omega$  be representing sequences consisting
  of projections such that 
 $e_n$ and $f_n$ are equivalent for any $n$. 
 Then there exist a representing sequence $(v_n)$ for 
 $V$ such that $v^*_nv_n=e_n$, $v_nv_n^*=f_n$. \\
 $(3)$ Let $\{E_{ij}\}_{1\leq i,j,\leq m}\subset M^\omega$ be a system of
  matrix units. Then there 
  exists a representing sequence $E_{ij}=\{e_{ij}^n\}$ such that
  $\{e_{ij}^n\}_{1\leq i,j\leq m}$ is a system of matrix units for every $n$.
 \end{lem}

 \section{Cohomology vanishing}\label{sec:coho}
 In this section, we mainly deal with actions of $\hG$ on factors of type II$_1$.
 However many parts of results in this section are valid for general
 factors (or von Neumann algebras). 

 We begin with 
 the following lemma, which is  known as the ``push-down lemma'' in subfactor
 theory \cite[Lemma 9.26]{EK-book}.
 \begin{lem}\label{lem:push}
  Let $\alpha$ be an action of $\hG$, and set
 $e:=1/|G|\sum_{\pi,i}d\pi\lambda_{{\pi_{ii}}}\in M\rtimes_\alpha \hG$. 
 For any $a\in M\rtimes_\alpha \hG$ 
  there exists $b\in M$ such that $ae=be$.
 \end{lem}
 \textbf{Proof.}  Let
 $a=\sum_{\rho,k,l}a_{\rho_{kl}}\lambda_{\rho_{kl}}$, $a_{\pi_{kl}}\in M$
 be an expansion of $a$. Then
 \begin{eqnarray*}
  |G|ae&=&\sum_{\pi,i,\rho,k,l}d\pi a_{\rho_{kl}}\lambda_{\rho_{kl}}\lambda_{\pi_{ii}} \\
 &=&\sum_{\pi,i,\rho,k,l}\sum_{\sigma,m,n,e}d\pi a_{\rho_{kl}}T_{\pi_i,\rho_k}^{\sigma_m,e}
 \lambda_{\sigma_{mn}}\overline{T_{\pi_i,\rho_l}^{\sigma_n,e}} \\
 &=&\sum_{\rho,k,l,\sigma,m,n}d\sigma a_{\rho_{kl}}\left(\sum_{\pi,i,e}
 \overline{T_{\sigma_m,\bar{\rho}_k}^{\pi_i,e}}
 T_{\sigma_n,\bar{\rho_{l}}}^{\pi_i,e}\right)\lambda_{\sigma_{mn}} \\
 &=&\sum_{\rho,k,\sigma,m}d\sigma a_{\rho_{kk}}\lambda_{\sigma_{mm}}\\
 &=&(\sum_{\rho,k}a_{\rho_{kk}})|G|e
 \end{eqnarray*}
 holds.
 Set $b:=\sum_{\rho,k}a_{\rho_{kk}}$, then we have $ae=be$ and $b\in M$. \hfill$\Box$

 \begin{prop}\label{prop:1vanish}
 Let $\alpha$ be an outer action of $\hG$. Then 
 any 1-cocycle for $\alpha$ is a coboundary. 
 \end{prop}
 \textbf{Proof.} Let $\{w_\pi\}$ be a 1-cocycle for $\alpha$, and
 $\lambda_\pi$ an implementing unitary in $M\rtimes_\alpha \hG$.
 It follows that $\{w_\pi\lambda_\pi\}$ is a representation of $\hG$.
 Set
 $e:=|G|^{-1}\sum_{\pi,i}d\pi\lambda_{\pi_{ii}}$, and
 $f:=|G|^{-1}\sum_{\pi,i}d\pi (w_{\pi}\lambda_{\pi})_{ii}=
 |G|^{-1}\sum_{\pi,i,j}d\pi w_{\pi_{ij}}\lambda_{\pi_{ji}}$. 
 Then $e$ and $f$ are projections in
 $M\rtimes_\alpha \hG$ with $E_M(e)=E_M(f)=|G|^{-1}$, where $E_M$ is the
 canonical conditional expectation 
 from $M\rtimes_\alpha \hG$ to $M$. 
 Since $M\rtimes_\alpha \hG$ is a factor due to the outerness of $\alpha$, 
 there exists $v\in M\rtimes_\alpha
 \hG$ such that $vev^*=f$. By Lemma \ref{lem:push}, we may assume $v\in M$.
  Since $vev^*=|G|^{-1}\sum_{\pi,i}d\pi v\lambda_{\pi_{ii}}v^*=
 |G|^{-1}\sum_{\pi,i,j}d\pi v\alpha_\pi(v^*)_{ij}\lambda_{\pi_{ji}}$, 
 we have $w_{\pi_{ij}}=v\alpha_{\pi}(v^*)_{ij}$, and this implies 
 $w_\pi=(v\otimes 1)\alpha_{\pi}(v^*)$. Especially $v$ is a unitary. \hfill$\Box$

 \begin{cor}\label{cor:equmu}
  Let $\alpha$ be an outer action of $\hG$ on $M$.
  Then there exists an $\alpha$-equivariant system of matrix units
  $\{e_{\pi_{ij},\rho_{kl}}\}\subset M$.
 \end{cor}
 \textbf{Proof.} We can choose  a system of matrix units $F=\{f_{\pi_{ij},\rho_{kl}}\}\subset
 M^\alpha$, since $M^\alpha$ is of type II$_1$. Then $\lambda^F_{\pi}$ is
 a 1-cocycle for $\alpha_\pi$. By Proposition \ref{prop:1vanish}, there
 exists $v\in U(M)$ such that $\lambda_\pi^F=(v^*\otimes
 1)\alpha_{\pi}(v)$. Define
 $E=\{e_{\pi_{ij},\rho_{kl}}\}:=\{vf_{\pi_{ij},\rho_{kl}}v^*\}$. 
 Then $\lambda_\pi^E=(v\otimes 1)\lambda_\pi^F(v^*\otimes 1)$ and 
 \begin{eqnarray*}
 \alpha_\pi(e_{\pi_{ij},\rho_{kl}})&=&\alpha_\pi(v)
 (f_{\pi_{ij},\rho_{kl}}\otimes 1)\alpha_{\pi}(v^*)\\
 &=&
 (v\otimes 1)\lambda_\pi^F(f_{\pi_{ij},\rho_{kl}}\otimes
 1)\lambda_\pi^{F*}(v^*\otimes 1) \\
 &=&\lambda_\pi^E(e_{\pi_{ij},\rho_{kl}}\otimes 1
 )\lambda_\pi^{E*}.
 \end{eqnarray*}
 \hfill$\Box$

 \begin{df}\label{df:2-cocycle} Let $\alpha_\pi\in \Mor(M, M\otimes B(H_\pi))$, 
 normalized as $\alpha_{\mathbf{1}}=\id_M$,
 and 
 $U_{\pi,\rho}\in U(M\otimes B(H_\pi)\otimes B(H_\rho))$.
 We say $\{\alpha_\pi, U_{\pi,\rho}\}_{\pi,\rho\in \hG}$ is a cocycle
  twisted 
  action of 
 $\hG$ if we have the 
  following.  \\
 $(1)$ $U_{\pi, \mathbf{1}}=U_{\mathbf{1},\pi}=1$. \\
 $(2)$ $\alpha_\pi\otimes \id_\rho\circ\alpha_\rho(x)U_{\pi,\rho}T=
 U_{\pi,\rho}T\alpha_\sigma(x), T\in (\sigma,\pi\otimes\rho)$. \\
 $(3)$ Let $\{T_{\pi,\rho}^{\sigma,a}\}_{a=1}^{N_{\pi\rho}^\sigma}$ be 
 an orthonormal basis for $(\sigma,\pi\otimes\rho)$. 
 Set $U_{\pi,\rho}^{\sigma,a}:=U_{\pi,\rho}T_{\pi,\rho}^{\sigma,a}$. 
 Then 
 $$
 \sum_{\eta,a,\xi,b}(\alpha_\pi\otimes \id)(U_{\rho,\sigma}^{\eta,a})
 U_{\pi,\eta}^{\xi,b} T^{\xi,b *}_{\pi,\eta}(1_\pi\otimes T_{\rho,\sigma}^{\eta,a*})=
 \sum_{\zeta,c,\xi,d}(U_{\pi,\rho}^{\zeta,c}\otimes 1_\sigma) U_{\zeta, \sigma}^{\xi,d} 
 T^{\xi,d*}_{\zeta,\sigma}(T_{\pi,\rho}^{\zeta,c*}\otimes 1_\sigma)$$
 holds. The unitary $U_{\pi,\rho}$ is called a 2-cocycle for $\alpha$.
 \end{df}

 We explain the meaning of Definition \ref{df:2-cocycle}(3). Note that 
 $\{(1_\pi \otimes T_{\rho,\sigma}^{\eta,a})T_{\pi,\eta}^{\xi,b}\}_{\eta,a,b}$ and 
 $\{(T_{\pi,\rho}^{\zeta,c}\otimes 1_\sigma)T_{\zeta,\sigma}^{\xi,d}\}_{\zeta,c,d}$
 are both orthonormal basis for $(\xi,\pi\otimes\rho\otimes \sigma)$. 
 Then $V_{(\zeta,c,d),(\eta,a,b)}=
 T_{\zeta,\sigma}^{\xi,d*}(T_{\pi,\rho}^{\zeta,c*}\otimes 1_\sigma)
 (1_\pi \otimes T_{\rho,\sigma}^{\eta,a})T_{\pi,\eta}^{\xi,b}
 \in (\xi,\xi)=
 \mathbf{C}$, and $V=\{V_{(\zeta,c,d),(\eta,a,b)}\}$ gives a unitary
 transformation between the above two orthonormal basis.
 From the condition (3), we get 
 $$(\alpha_\pi\otimes \id)(U_{\rho,\sigma}^{\eta,a})
 U_{\pi,\eta}^{\xi,b} 
 =\sum_{\zeta,c,d}
 (U_{\pi,\rho}^{\zeta,c}\otimes 1_\sigma) U_{\zeta, \sigma}^{\xi,d} 
 V_{(\zeta,c,d),(\eta,{a,b})}.$$
 This shows  that the same unitary $V$ gives 
 the transformation between $(\alpha_\pi\otimes \id)(U_{\rho,\sigma}^{\eta,a})
 U_{\pi,\eta}^{\xi,b} $ and 
 $(U_{\pi,\rho}^{\zeta,c}\otimes 1_\sigma) U_{\zeta, \sigma}^{\xi,d} $. 

 Let $\{\alpha_\pi, U_{\pi,\sigma}\}$ be a cocycle twisted action of
 $\hG$, and extend to that of $\rG$ as in the remark in \S\ref{sub:inter}. Then 
 we have $(\alpha_\pi\otimes \id_\rho)\circ \alpha_\rho=\Ad
 (U_{\pi,\rho})\alpha_{\pi\otimes \rho}$, and 
 $(U_{\pi,\rho}\otimes 1_\sigma)U_{\pi\otimes \rho,\sigma}=
 (\alpha_\pi\otimes \id_\rho\otimes \id_\sigma)
 (U_{\rho,\sigma})U_{\pi,\rho\otimes \sigma}$ from the above equalities.

 \begin{df}\label{df:2coboundary}
  Let $\{\alpha_\pi, U_{\pi,\rho}\}$ be a cocycle twisted action of $\hG$. 
 We say that $U_{\pi,\rho}$ is a coboundary if there exist unitaries $W_\pi\in M\otimes 
 B(H_\pi)$, normalized as $W_\mathbf{1}=1$,
 such that $$W_\pi\alpha_\pi\otimes \id(W_\rho)U_{\pi,\rho}T=
 T W_\sigma, T\in (\sigma,\pi\otimes\rho).$$
 Define $$(\partial_\alpha W)_{\pi,\rho}:=
 \alpha_\pi\otimes \id_\rho(W_\rho)
 (W_\pi\otimes 1_\rho)W_{\pi\otimes \rho}^*=
 \sum_{\sigma,e}\alpha_\pi\otimes \id_\rho(W_\rho)
 (W_\pi\otimes 1_\rho)T_{\pi,\rho}^{\sigma,e} W^*_\sigma T_{\pi,\rho}^{\sigma,e*}.$$
 Then the above condition is shown to be  
 equivalent to $U_{\pi,\rho}=(\partial_\alpha W^*)_{\pi,\rho}$.
 \end{df}

 Let $\{\alpha, U_{\pi,\rho}\}$ be a cocycle twisted action, and assume
  $U_{\pi,\rho}=(\partial_\alpha W^*)_{\pi,\rho}$ for some $\{W_\pi\}$. 
  Then $\Ad W_\pi\alpha_\pi$
  becomes a genuine action of $\hG$. \\

 \noindent
 \textbf{Remark.} Here we make a useful remark on perturbation of cocycle
 twisted actions.
 Let $\alpha$ be an action of $\hG$, and $w_\pi\in U(M\otimes
 B(H_\pi))$. Then $\tilde{\alpha}_\pi=\Ad w_\pi \alpha_\pi$ is a cocycle twisted
 action with a 2-cocycle $u(\pi,\rho)=\partial_{\tilde{\alpha}} w(\pi,\rho)$.
 If there exists another unitary $\bar{w}_\pi$ such that
 $\partial_{\tilde{\alpha}}\bar{w}^*(\pi,\rho )=u(\pi,\rho)$, then
 it is easy to verify that $\bar{w}_\pi w_\pi$ is a 1-cocycle for
 $\alpha$. \\

 In a similar way as in \cite{J-act} and \cite{Su-homoII}, we can prove the 2-cohomology
 vanishing theorem for cocycle twisted actions of $\hG$ as follows.

 \begin{thm}\label{thm:2vanish-1}
 Let $\{\alpha_\pi, U_{\pi,\rho}\}$ be a $($not necessary outer$)$ 
 cocycle twisted action of $\hG$. Then $U_{\pi,\rho}$ 
 is a coboundary. 
 \end{thm}
 \textbf{Proof.} 
 Fix a finite dimensional subfactor $K\subset M$,
 $K\cong M_{|G|}(\mathbf{C})$, and a system of matrix units $\{e_{ij}\}_{1\leq
 i,j\leq |G|}$ for $K$. Choose a unitary $u_\pi\in M\otimes B(H_\pi)$
 with $\Ad u_\pi\alpha_\pi(e_{ij})=e_{ij}\otimes 1$, and
 set $\tilde{\alpha}_\pi:=\Ad u_\pi \alpha_\pi$. 
 Then $\tilde{\alpha}_\pi(x)=x\otimes 1$, $x\in K$. Hence
 $\tilde{\alpha}_\pi$ sends $K'\cap M$ into $(K'\cap M)\otimes
 B(H_\pi)$. Moreover if we define
 $\tilde{U}_{\pi,\rho}:=u_\pi^{12}\alpha_\pi(u_\rho)U_{\pi,\rho}u_{\pi\otimes
 \rho}^*=\sum_{\sigma,a} u_\pi^{12}\alpha_\pi(u_\rho)U_{\pi,\rho}^{\sigma,a}
 u_\sigma^* T_{\pi\rho}^{\sigma,a*}$, then $\tilde{U}_{\pi,\rho}\in (K'\cap M)\otimes
 B(H_\pi)\otimes B(H_\rho)$ and $\{\tilde{\alpha},\tilde{U}_{\pi,\rho}\}$
 is a cocycle twisted action of $\hG$ on $K'\cap M$. 
 It is trivial that  
 $\tilde{U}_{\pi,\rho}$ is a coboundary if and only if  so is $U_{\pi,\rho}$. 

 Hence  we may assume $\alpha$ is of the form
 $\alpha_\pi=\alpha_\pi^0\otimes \id$ on $N\otimes B(\ell^2(\hG))$ and
 $U_{\pi,\rho}=\sum_{\pi,i,j,\rho,k,l} 
 U_{\pi_{ij},\rho_{kl}}\otimes 1 \otimes e^\pi_{ij}\otimes e^\rho_{kl} 
 \in N\otimes \mathbf{C}1\otimes B(H_\pi)\otimes
 B(H_\rho)$. Fix a system of matrix units $\{f_{\xi_{ab},\eta_{cd}}\}$ for
 $B(\ell^2(\hG))$. 

 In the rest of this section, we denote
 $U_{\pi,\rho}^{\sigma,a}$ and  
 $T_{\pi,\rho}^{\sigma,a}$ by $U_{\pi,\rho}^{\sigma}$ and
 $T_{\pi,\rho}^{\sigma}$ respectively to simplify notations.

 Define $w_{\pi_{ij}}$ as 
 $$w_{\pi_{ij}}:=\sum_{\xi,a,b,\eta,c,d}\sqrt{\frac{d\xi}{d\eta}}T_{\pi_i\xi_a}^{\eta_c}
 U_{\pi_j\xi_b}^{\eta_d^*}\otimes f_{\eta_{cd},\xi_{ab}}\in N\otimes B(\ell^2(\hG)).$$
 We will see $w_\pi\in U(M\otimes B(H_\pi))$, and $(\partial_\alpha
 w^*)_{\pi,\rho}=U_{\pi,\rho}$. 
 At first we verify that $w_\pi$ is a unitary. We will see  
 $\sum_k w_{\pi_{ik}}w_{\pi_{jk}}^*=\delta_{i,j}$ holds as follows. 

 \begin{eqnarray*}
 \lefteqn{ 
 \sum_{k}w_{\pi_{ik}}w_{\pi_{jk}}^* }\\
 &=&
  \sum_k\left(\sum_{\rho,l,m, \sigma,s,t}\sqrt{\frac{d\rho}{d\sigma}}T_{\pi_i\rho_l}^{\sigma_s}
 U_{\pi_k\rho_m}^{\sigma_t^*}\otimes f_{\sigma_{st},\rho_{lm}}\right)
  \left(\sum_{\xi,a,b,\eta,c,d}\sqrt{\frac{d\xi}{d\eta}}\overline{T_{\pi_j\xi_a}^{\eta_c}}
 U_{\pi_k\xi_b}^{\eta_d}\otimes f_{\xi_{ab},\eta_{cd}}\right) \\
 &=&
 \sum_{\rho,l, \sigma,s,t,\eta,c,d}\frac{d\rho}{\sqrt{d\sigma d\eta}}
 T_{\pi_i\rho_l}^{\sigma_s}
 \overline{T_{\pi_j\rho_l}^{\eta_c}}
 \left(\sum_{k,m}U_{\pi_k\rho_m}^{\sigma_t^*}
 U_{\pi_k\rho_m}^{\eta_d}\right)\otimes f_{\sigma_{st},\eta_{cd}} \\
 &=&
 \sum_{\rho,m, \eta,c,d,s}\frac{d\rho}{d\eta}T_{\pi_i\rho_l}^{\eta_s}
 \overline{T_{\pi_j\rho_l}^{\eta_c}}
 \otimes f_{\eta_{sd},\eta_{cd}} \\
 &=&
 \sum_{\eta,c,d,s}
 \left(\sum_{\rho,l}\overline{T_{\bar{\pi}_i\eta_s}^{\rho_l}}
 T_{\bar{\pi}_j\eta_c}^{\rho_l}\right)
 \otimes f_{\eta_{sd},\eta_{cd}} \\
 &=&\delta_{i,j}
 \sum_{\eta,c,d}
 1 \otimes f_{\eta_{cd},\eta_{cd}} \\
 &=&\delta_{i,j}.
 \end{eqnarray*}
 In a similar way as above,
 $\sum_{k}w_{\pi_{ki}}^*w_{\pi_{kj}}=\delta_{i,j}$ can be verified. Hence
 $w_\pi$ is indeed a unitary.

 We next show that $U_{\pi,\rho}=(\partial_\alpha w^*)_{\pi,\rho}$.
 It suffices to show
 $w_\pi\alpha_{\pi}(\omega_\rho)U_{\pi,\rho}^{\sigma}
 w_\sigma^*=T_{\pi,\rho}^{\sigma}$. 
 This follows from  the computation.
 \begin{eqnarray*}
  \lefteqn{(w_\pi\alpha_{\pi}(\omega_\rho)U_{\pi,\rho}^\sigma 
 w_\sigma^*)_{\pi_i,\rho_k}^{\sigma_m}} \\ 
 &=& \sum_{j,l,n,a} w_{\pi_{ij}}\alpha_{\pi}(w_{\rho_{kl}})_{jn}
 (U_{\pi_{n}\rho_{l}}^{\sigma_{a}}\otimes 1)
 w_{\sigma_{ma}}^* \\
 &=&\sum_{j,l,n,a,\xi,b,c,\atop \eta,d, e,\zeta,s,t, \phi,u,v}
 \sqrt{\frac{d\xi}{d\eta}} \sqrt{\frac{d\zeta}{d\xi}}\sqrt{\frac{d\zeta}{d\phi}}
 T_{\pi_{i}\xi_b}^{\eta_d}U_{\pi_j\xi_{c}}^{\eta_e*}
 T_{\rho_k\zeta_s}^{\xi_{b}}
  \alpha_{\pi}(U_{\rho_{l}\zeta_t}^{\xi_c*})_{jn}
 U_{\pi_{n}\rho_{l}}^{\sigma_{a}}
 \overline{T_{\sigma_m, \zeta_s}^{\phi_u}}U_{\sigma_{a}\zeta_t}^{\phi_{v}}
 \otimes f_{\eta_{ed},\phi_{uv}} \\
 &=&\sum_{l,n,a,\eta,d, e,\atop \zeta,s,t,\phi,u,v}
 \frac{d\zeta}{\sqrt{d\eta d\phi}}
 \left(\sum_{j,\xi,b,c}
  \alpha_{\pi}(U_{\rho_{l}\zeta_t}^{\xi_c})_{nj}U_{\pi_j\xi_{c}}^{\eta_e}
 \overline{T_{\rho_k\zeta_s}^{\xi_{b}}}
 \overline{T_{\pi_{i}\xi_b}^{\eta_d}}
 \right)^*
 U_{\pi_{n}\rho_{l}}^{\sigma_{a}}
 \overline{T_{\sigma_m, \zeta_s}^{\phi_u}}U_{\sigma_{a}\zeta_t}^{\phi_{v}}
 \otimes f_{\eta_{ed},\phi_{uv}} \\
 &=&\sum_{l,n,a,\eta,d, e,\atop \zeta,s,t,\phi,u,v}
 \frac{d\zeta}{\sqrt{d\eta d\phi}}
 \left(\sum_{\xi,b,c}
 U_{\pi_n\rho_l}^{\xi_c}U_{\xi_{c}\zeta_t}^{\eta_e}
 \overline{T_{\pi_{i}\rho_k}^{\xi_b}}
 \overline{T_{\xi_b\zeta_s}^{\eta_d}}
 \right)^*
 U_{\pi_{n}\rho_{l}}^{\sigma_{a}}
 \overline{T_{\sigma_m, \zeta_s}^{\phi_u}}U_{\sigma_{a}\zeta_t}^{\phi_{v}}
 \otimes f_{\eta_{ed},\phi_{uv}} \\
 &=&\sum_{a,\xi,b,c,\eta,d, e,\atop \zeta,s,t,\phi,u,v}
 \frac{d\zeta}{\sqrt{d\eta d\phi}}
 {T_{\pi_{i}\rho_k}^{\xi_b}}
 T_{\xi_b\zeta_s}^{\eta_d}
 U_{\xi_{c}\zeta_t}^{\eta_e*}\left(\sum_{n,l}
 U_{\pi_n\rho_l}^{\xi_c*}
 U_{\pi_{n}\rho_{l}}^{\sigma_{a}}\right)
 \overline{T_{\sigma_m, \zeta_s}^{\phi_u}}U_{\sigma_{a}\zeta_t}^{\phi_{v}}
 \otimes f_{\eta_{ed},\phi_{uv}} \\
 &=&\sum_{b ,\eta,d, e,\atop \zeta,s, \phi,u,v}
 \frac{d\zeta}{\sqrt{d\eta d\phi}}
 {T_{\pi_{i}\rho_k}^{\sigma_b}}
 T_{\sigma_b\zeta_s}^{\eta_d}
 \overline{T_{\sigma_m, \zeta_s}^{\phi_u}}
 \sum_{a,t}\left(
 U_{\sigma_a\zeta_t}^{\eta_e*}
 U_{\sigma_{a}\zeta_t}^{\phi_{v}}\right)
 \otimes f_{\eta_{ed},\phi_{uv}} \\
 &=&\sum_{b ,\eta,d, e,u,\zeta,s}
 \frac{d\zeta}{d\eta}
 {T_{\pi_{i}\rho_k}^{\sigma_b}}
 T_{\sigma_b\zeta_s}^{\eta_d}
 \overline{T_{\sigma_m, \zeta_s}^{\eta_u}}
 \otimes f_{\eta_{ed},\eta_{ud}} \\
 &=&\sum_{b ,\eta,d, e,u}
 T_{\pi_{i}\rho_k}^{\sigma_b}
 \left(\sum_{\zeta,s}\overline{T_{\bar{\sigma}_b\eta_d}^{\zeta_s}}
 {T_{\sigma_m, \eta_u}^{\zeta_s}}\right)
 \otimes f_{\eta_{ed},\eta_{ud}} \\
 &=&
 T_{\pi_{i}\rho_k}^{\sigma_m}
 \sum_{\eta,d, e}
 1\otimes f_{\eta_{ed},\eta_{ed}} \\
 &=&
 T_{\pi_{i}\rho_k}^{\sigma_m}.
 \end{eqnarray*}

 Hence we have
 $U_{\pi,\rho}=(\partial_\alpha w^*)_{\pi,\rho}$. \hfill$\Box$

 We need
 another type of 2-cohomology vanishing theorem, which asserts that 
 we can choose a coboundary close to 1 if a 2-cocycle is close to 1. 

  \begin{thm}\label{thm:2vanish-2}
 Let $\{\alpha_\pi, U_{\pi,\rho}\}$ be a cocycle twisted outer action of
   $\hG$. If  $\|U_{\pi,\rho}-1\|_2<\delta$ for 
   sufficiently small enough $\delta$,  
 then there exist a unitary $W_\pi\in M\otimes B(H_\pi)$ such that 
 $U_{\pi,\rho}=(\partial_\alpha W^*)_{\pi,\rho}$ and
 $\|W_\pi-1\|_2<f(\delta)$.
 Here $f(\delta)$ is a
 positive valued function, which depends only on $G$ and is
   independent from $\alpha$ and $U_{\pi,\rho}$, 
   with $\lim_{\delta\rightarrow 0}f(\delta)=0$.
 \end{thm}
 \textbf{Proof.} Let $N:=M\rtimes_{\alpha, U}\hG$ be a twisted crossed
 product, and $\lambda_{\pi}$ be an implementing unitary. (See Appendix
 for the twisted crossed product construction.)
 Hence we have $\Ad \lambda_\pi(x\otimes 1_\pi)=\alpha_\pi(x)$ for $x\in M$, 
 $\lambda_\pi^{12}\lambda_\rho^{13}T_{\pi,\rho}^\sigma=U_{\pi,\rho}T_{\pi,\rho}^{\sigma}
 \lambda_\sigma$ for $T_{\pi,\rho}^\sigma\in(\sigma,\pi\otimes \rho)$,
 and $M\rtimes_{\alpha,U}\hG=M\vee\{\lambda_{\pi_{ij}}\}$. 
 By Theorem \ref{thm:2vanish-1}, there exist a unitary $w_\pi\in M\otimes
 B(H_\pi)$, $\pi\in\hG$ such 
 that $U_{\pi,\rho}=(\partial_\alpha w^*)_{\pi,\rho}$.
 This implies
 $\tilde{\lambda}_{\pi}:=w_{\pi}\lambda_{\pi}$
 is a representation of $\hG$. Thus $e:=|G|^{-1}\sum_{\pi,i}d\pi
 \tilde{\lambda}_{\pi_{ii}}$ is a projection with $E_M(e)=1/|G|$.

 (If one is not
 familiar to the twisted crossed product, he (or she) may treat it as
 follows. Let $w_\pi$ be as above.
 Since $\Ad w_\pi\alpha_\pi$ is a usual action, we can construct
 a usual crossed product algebra $M\rtimes_{\Ad w\alpha}\hG$. Let
 $\tilde{\lambda}_\pi$ be an implementing unitary, and set
 $\lambda_\pi:=w_\pi^*\lambda_\pi$. Then it is easy to see $\{\lambda_\pi\}$
 behave like as the implementing unitary in $M\rtimes_{\alpha,U}\hG$.
 Hence $M\vee \{\lambda_\pi\}$ is
 identified with the twisted crossed product $M\rtimes_{\alpha, U}\hG$. )

 Set $f:=|G|^{-1}\sum_{\pi,i}d\pi\lambda_{\pi_{ii}}$. 
 We will show $f$ is
 almost a projection, and apply Lemma \ref{lem:nearproj}. 
 Set $\Lambda_\pi:=\sum_i\lambda_{\pi_{ii}}$ and we
 investigate
 $\Lambda_\pi\Lambda_\rho-\sum_\sigma N_{\pi,\rho}^\sigma\Lambda_\sigma$
 at first.

 Since we have
 \begin{eqnarray*}
  \Lambda_\pi\Lambda_\rho&=&\sum_{i,k}\lambda_{\pi_{ii}}\lambda_{\rho_{kk}}\\
 &=&\sum_{i,k,\sigma,m,n,a}U_{\pi_i,\rho_k}^{\sigma_m,a}\lambda_{\sigma_{mn}}
 \overline{T_{\pi_i,\rho_k}^{\sigma_n,a}} \\
 &=&
 \sum_{i,j,k,l,\sigma,m,n}U_{\pi_{ij},\rho_{kl}}T_{\pi_j,\rho_l}^{\sigma_m,a}\lambda_{\sigma_{mn}}
 \overline{T_{\pi_i,\rho_k}^{\sigma_n,a}} \\
 &=&
 \sum_{i,j,k,l,\sigma,m,n}(U_{\pi_{ij},\rho_{kl}}-\delta_{ij}\delta_{kl})
 T_{\pi_j,\rho_l}^{\sigma_m,a}\lambda_{\sigma_{mn}}
 \overline{T_{\pi_i,\rho_k}^{\sigma_n,a}} +
 \sum_{i,k,\sigma,m,n}T_{\pi_i,\rho_k}^{\sigma_m,a}\lambda_{\sigma_{mn}}
 \overline{T_{\pi_i,\rho_k}^{\sigma_n,a}} \\
 &=&
 \sum_{i,j,k,l,\sigma,m,n}(U_{\pi_{ij},\rho_{kl}}-\delta_{ij}\delta_{kl})
 T_{\pi_j,\rho_l}^{\sigma_m}\lambda_{\sigma_{mn}}
 \overline{T_{\pi_i,\rho_k}^{\sigma_n}} +
 \sum_\sigma N_{\pi,\rho}^\sigma \Lambda_\sigma,
 \end{eqnarray*}
 we get the following estimate.
 \begin{eqnarray*}
  \left\|\Lambda_\pi\Lambda_\rho - \sum_\sigma N_{\pi,\rho}^\sigma\Lambda_\sigma\right\|_2&\leq &
 \sum_{i,j,k,l}\left\|(U_{\pi_{ij}\rho_{kl}}-\delta_{ij}\delta_{kl})
 \sum_{\sigma,m,n,a}T_{\pi_j\rho_l}^{\sigma_m,a}\lambda_{\sigma_{mn}}
 \overline{T_{\pi_i\rho_k}^{\sigma_n,a}}\right\|_2 \\
 &\leq &
 \sum_{i,j,k,l}\|U_{\pi_{ij}\rho_{kl}}-\delta_{ij}\delta_{kl}\|_2 \\
 &\leq & d\pi d\rho
 \sqrt{\sum_{i,j,k,l}\|U_{\pi_{ij}\rho_{kl}}-\delta_{ij}\delta_{kl}\|_2^2 }\\
 &\leq & d\pi d\rho\sqrt{d\pi d\rho}\delta.
 \end{eqnarray*}

 Now we give the estimate of $\|f^2-f\|_2$. Since
 \begin{eqnarray*}
  f^2-f&=&\frac{1}{|G|^2}\sum_{\pi,\rho}d\pi d\rho\Lambda_\pi\Lambda_\rho-
 \frac{1}{|G|}\sum_\sigma d\sigma\Lambda_\sigma \\
 &=&
 \frac{1}{|G|^2}\sum_{\pi,\rho,\sigma}d\pi d\rho\left(\Lambda_\pi\Lambda_\rho-
 N_{\pi,\rho}^\sigma\Lambda_\sigma\right)+
 \frac{1}{|G|^2}\sum_{\pi,\rho,\sigma}d\pi d\rho N_{\pi,\rho}^\sigma\Lambda_\sigma-
 \frac{1}{|G|}\sum_\sigma d\sigma\Lambda_\sigma \\
 &=&
 \frac{1}{|G|^2}\sum_{\pi,\rho,\sigma}d\pi d\rho\left(\Lambda_\pi\Lambda_\rho-
 N_{\pi,\rho}^\sigma\Lambda_\sigma\right)
 \end{eqnarray*}
 holds, we get
 \begin{eqnarray*}
 \|f^2-f\|_2&\leq &\frac{1}{|G|^2}\sum_{\pi,\rho}d\pi d\rho
 \left\|\Lambda_\pi\Lambda_\rho-\sum_\sigma N_{\pi,\rho}^\sigma\Lambda_\sigma \right\|_2 \\ 
 &\leq &\frac{1}{|G|^2}\sum_{\pi,\rho}d\pi^2d\rho^2\sqrt{d\pi d\rho}\delta.
 \end{eqnarray*}

 Next we estimate $\|f^*-f\|_2$. 
 Set $\tilde{U}_{\pi_{i},\bar{\pi}_j}=\sum_k
 U_{\pi_{ik},\bar{\pi}_{jk}}$.
 Then
 $\lambda_{\bar{\pi}_{ij}}^*=\sum_k\tilde{U}_{\pi_{k},\bar{\pi}_i}^*\lambda_{\pi_{kj}}$.
 (See Appendix), and we have
 $f^*=|G|^{-1}\sum_{\pi,i,j}d\pi\tilde{U}_{\pi_{j},\bar{\pi}_i}^*\lambda_{\pi_{ji}}$. 

 Hence 
 \begin{eqnarray*}
 \|f^*-f\|_2&\leq 
  &\frac{1}{|G|}\sum_{\pi,i,j}d\pi\left\|(\tilde{U}_{\pi_{j},\bar{\pi}_i}^*-\delta_{i,j})
 \lambda_{\pi_{ji}}\right\|_2\\ 
 &\leq  &
 \frac{1}{|G|}\sum_{\pi,i,j,k}d\pi\left\|(U_{\pi_{jk},\bar{\pi}_{ik}}^*-
 \delta_{i,k}\delta_{j,k})\right\|_2\\ 
 &\leq&
 \frac{1}{|G|}\sum_\pi d\pi^2\sqrt{d\pi}\sqrt{\sum_{i,j,k}\left\|(U_{\pi_{jk},\bar{\pi}_{ik}}^*-
 \delta_{i,k}\delta_{j,k})\right\|_2^2}\\ 
 &\leq&
 \frac{1}{|G|}\sum_\pi d\pi^3\sqrt{d\pi}\delta
 \end{eqnarray*}
 holds.

 If we set $C:=\max\{|G|^{-1}\sum_{\pi,\rho}d\pi^2d\rho^2\sqrt{d\pi
 d\rho}, |G|^{-1}\sum_\pi d\pi^3\sqrt{d\pi}\}$, then 
 $\|f^2-f\|_2\leq C\delta$ and $\|f^*-f\|_2\leq C\delta$.
 Here note that $C$ is determined only by $G$.

 By Lemma \ref{lem:nearproj}, there exists a projection $p\in
 M\rtimes_{\alpha, U}\hG$ such that 
 $\|f-p\|_2<6\sqrt[4]{C\delta}$ and $\tau(p)=1/|G|$, provided $C\delta <1/4$. Then 
 there exists a unitary $a\in M\rtimes_{\alpha,U}\hG$ such that $p=aea^*$, and by
 Lemma \ref{lem:push}, there exists $u \in M$ with $p=ueu^*$.
 Note that $u$ is not necessary
 a unitary, and we have $\|u\|\leq |G|$ since $u=|G|E_M(ae)$.
 Then
 $\|f-ueu^*\|_2<6\sqrt[4]{C\delta}$. By applying the canonical
 conditional expectation, we get
 $\|1-uu^*\|_2<6|G|\sqrt[4]{C\delta}$. By Lemma \ref{lem:nearunit}, there exists a
 unitary $v\in M$ such that $\|u-v\|_2<6|G|(3+|G|)\sqrt[4]{C\delta}$. 
 Hence we have $\|f-vev^*\|_2<f(\delta)$ for some positive valued
 function $f(\delta)$, which depends only on $G$, with
 $\lim_{\delta\rightarrow 0}f(\delta)=0$. 
 Set $\bar{w}_\pi:=(v\otimes
 1)w_\pi\alpha_\pi(v^*)$. Then $\bar{w}_\pi$ is a coboundary for
 $U_{\pi,\rho}$, and $\|\bar{w}_\pi-1\|_2<f(\delta)$ holds by looking
 at coefficients of
 $f-vev^*=|G|^{-1}\sum_{\pi,i,j}d\pi(1-\bar{w}_{\pi_{ij}})
 \lambda_{\pi_{ji}}$.
 \hfill$\Box$

 \section{Actions and ultra product}\label{sec:ultra}
 Fix a free ultrafilter $\omega$ over $\mathbf{N}$. 
 Then $\alpha_\pi^\omega$ is an action of $\hG$ on $M^\omega$, however
 $\alpha_\pi^\omega $ does not preserve $M_\omega$ unless $\hG$ is a
 group, i.e., $G$ is commutative.

 Our first task in this section is to
  show the existence of an outer action of
 $\hG\times \hG$ on $M_\omega$ by modifying $\alpha_\pi^\omega$.

 We first recall  Ocneanu's central freedom lemma.
 \begin{lem}\label{lem:centfree}
  Let $A\subset B\subset C$ be finite von Neumann algebras with $A\cong \cR$.
  Then $(A'\cap B^\omega)'\cap C^\omega=A\vee (B'\cap C)^\omega$.
 \end{lem}
 See \cite{EK-book} for proof.

 From now on, we assume $M\cong\cR$.  
 Set $M_1:=M\rtimes_\alpha \hG$. Let $\lambda_\pi\in M_1\otimes B(H_\pi)$
 be an implementing unitary for $\alpha$.  

 \begin{lem}\label{lem:approinner}
  Let $H$ be a finite dimensional Hilbert space.
 Then $\alpha \in \Mor(\cR, \cR\otimes B(H))$ is approximately 
  inner in the following sense; there exists a sequence of unitary
  $\{u_n\}\subset \cR\otimes B(H)$ such that 
 $\lim_n\|\Ad u_n(x\otimes 1)-\alpha(x)\|_{2}=0$, $x\in \cR$.
  \end{lem}
 \textbf{Proof.} Represent $\cR=\bigotimes_{n=1}^\infty M_2(\mathbf{C})$,
 and set $L_n:=\bigotimes_{k=1}^nM_2(\mathbf{C})$.
 Let $\{e_{ij}^n\}$ be a system of
 matrix units for $L_n$. Then $\{e_{ij}^n\otimes 1\}$ and $\{\alpha(e_{ij}^n)\}$
 are both systems of matrix units in $\cR\otimes B(H)$. Hence there exists
 a unitary $u_n\in \cR\otimes B(H)$ with $\alpha(e_{ij})=\Ad
 u_n(e_{ij}\otimes 1)$, and hence $\alpha(x)=\Ad u_n(x\otimes 1)$ for
 $x\in L_n$. Then it is easy to see $\lim_n\|\alpha(x)-\Ad u_n(x\otimes
 1)\|_2=0$ 
 for $x\in \cR$. \hfill$\Box$ 

 By Lemma \ref{lem:approinner}, there exists a unitary $U_\pi\in M^\omega\otimes B(H_\pi)$
 such that
 $\alpha_\pi(x)=\Ad U_\pi(x\otimes 1)$ for $x\in M\subset M^\omega$,
 $\pi\in \hG$.
 Set $V_\pi:=U_\pi^*\lambda_\pi\in M_1^\omega$.

 \begin{lem}\label{lem:cocycleaction}
   Define $\gamma^1_\pi(x):=\Ad V_\pi(x\otimes  1)$, 
 $\gamma^2_\rho(x):=\Ad U_\rho^*(x\otimes 1_\rho)$. 
  Then $\gamma_{\pi\hat{\otimes}\rho}:=(\gamma^1_\pi\otimes 1_\rho)\circ\gamma^2_\rho$ 
 defines an outer cocycle twisted
  action of $\hG\times \hG$ on $M_\omega$. 
 \end{lem}
 \textbf{Proof.} Define $W_{\pi\hat{\otimes }\rho}:=
 V_\pi^{12}U^{13*}_\rho\in M_1^\omega\otimes B(H_\pi)\otimes B(H_\rho)$.
 Then
 $\gamma$ is a perturbation of the trivial action of $\hG\times \hG$ on $M_1^\omega$ by
 $W_{\pi\hat{\otimes }\rho}$, 
 and $w(\pi^1\hat{\otimes }\rho^1,
 \pi^2\hat{\otimes }\rho^2):=\partial_\gamma (W)(
 \pi^1\hat{\otimes }\rho^1,\pi^2\hat{\otimes }\rho^2)$
 is a 2-cocycle for $\gamma$. Hence we only have to verify 
 $\gamma$ preserves $M_\omega$,
 outer on $M_\omega$, and 
 $w(\pi^1\hat{\otimes }\rho^1,\pi^2\hat{\otimes }\rho^2)\in 
 M_\omega\otimes
 B(H_{\pi^1\hat{\otimes} \rho^1})\otimes B(H_{\pi^2\hat{\otimes}\rho^2}) $. 

 We verify that $\gamma^i_\pi\in\Mor(M_\omega, 
 M_\omega\otimes B(H_\pi))$, $i=1,2$.
 Indeed this  follows from the computation below, 
 where $x\in M, a\in M_\omega$. Note $\gamma_\pi^1(a)=\Ad
 U_\pi^*\alpha_\pi^\omega(a)$ for $a\in M_\omega$.
 $$U_\pi^*(a\otimes 1)U_\pi(x\otimes 1)=U_\pi^*(a\otimes 1)\alpha_\pi(x)U_\pi
 =U_\pi^*\alpha_\pi(x)(a\otimes 1)U_\pi=(x\otimes 1)\Ad U_\pi^*(a\otimes 1),$$
 $$U_\pi^*\alpha_\pi^\omega(a)U_\pi(x\otimes 1)=
 U_\pi^*\alpha_\pi^\omega(a)\alpha_\pi(x)U_\pi=
 U_\pi^*\alpha_\pi^\omega(ax)U_\pi=
 U_\pi^*\alpha_\pi^\omega(xa)U_\pi=(x\otimes 1)U_\pi^*\alpha_\pi^\omega(a)U_\pi.$$
 Then it is trivial $\gamma_{\pi\hat{\otimes}\rho}\in
 \Mor(M_\omega,M_\omega\otimes B(H_\pi)\otimes B(H_\rho))$.

 We next verify $\gamma_{\pi\hat{\otimes }\rho}$ is outer on $M_\omega$.
 We divide to $\pi\ne \mathbf{1} $ case and $\pi=\mathbf{1}$ case.
 Fix $\mathbf{1}\ne \pi\in \hG$ and assume
 $\gamma_{\pi\hat{\otimes}\rho}(x)a=a(x\otimes 1)$, $x\in M_\omega$,  
 holds for some $a\in M_\omega \otimes B(H_\pi)\otimes B(H_\rho)$. 
 On one hand, $b:=U_\rho^{13}V_\pi^{12*}a\in
  (M_\omega\otimes \mathbf{C})'\cap
 M_1^\omega\otimes B(H_\pi)\otimes B(H_\rho)=M\otimes B(H_\pi)\otimes
 B(H_\rho)$ 
 by Lemma \ref{lem:centfree}.
 On the other hand, it is easy to see $b=U_\rho^{13}V_\pi^{12*}a $ is of
 the form $\sum_{i,j,k,l,m,n}X^{k,l}_{i,j,m,n}\lambda_{\bar{\pi}_{m,n}}\otimes
 e^\pi_{ij}\otimes e^\rho_{kl}$, $X^{k,l}_{i,j,m,n}\in M^\omega$. Hence
 if $\pi\ne \mathbf{1}$, $b$ must be $0$, and consequently $a=0$.

 Assume $\pi=\mathbf{1}$, and we verify $\gamma_{\mathbf{1}\hat{\otimes
 }\rho}=\gamma^2_\rho$ is outer for $\rho\ne \mathbf{1}$. 
 Assume $\gamma^2_\rho(x) a=a(x\otimes 1)$. 
 Then  $b_\rho:=U_\rho a\in (M_\omega\otimes \mathbf{C})'\cap
 M^\omega\otimes B(H_\rho)=M\otimes B(H_\rho)$ by Lemma \ref{lem:centfree}.

 Then we have 
 $$\alpha_\rho(x)b_\rho=U_\rho(x\otimes 1)U_\rho^*b_\rho=U_\rho(x\otimes
 1)a=U_\rho a(x\otimes 1)=b_\rho(x\otimes 1)$$ 
 for $x\in M$. Since $\alpha$ is an outer action on $M$,  $b_\rho$ is $0$ and hence so
 is $a$.

 We will see $w(\pi^1\hat{\otimes}\rho^1,\pi^2\hat{\otimes}\rho^2)\in
 U(M_\omega\otimes 
 B(H_{\pi^1\hat{\otimes} \rho^1})\otimes B(H_{\pi^2\hat{\otimes}\rho^2}) )$. 
 Although we can prove this directly, 
 we will show the statement for
 $w(\pi\hat{\otimes}\mathbf{1},\rho\hat{\otimes}\mathbf{1})$, 
 $w(\mathbf{1}\hat{\otimes}\pi,\mathbf{1}\hat{\otimes}\rho)$ and  
 $w(\mathbf{1}\hat{\otimes}\pi, \rho\hat{\otimes}\mathbf{1})$ separately
 to abuse notations.
 Then we can obtain the desired result since 
 \begin{eqnarray*}
 \lefteqn{w(\pi^1\hat{\otimes}\rho^1,\pi^2\hat{\otimes}\rho^2)
 =\gamma_{{\pi_1}}^1\otimes \id_{\rho_1}\otimes \id_{\pi_2}\otimes \id_{\rho_2}
 (w(\mathbf{1}\hat{\otimes }\rho_1,\pi_2\hat{\otimes}\mathbf{1})\otimes
 1_{\rho_2})\times} \\
 &&(w(\pi_1\hat{\otimes}\mathbf{1},\pi_2\hat{\otimes}\mathbf{1})\otimes
 1_{\rho_1}\otimes 1_{\rho_2})
 \gamma_{\pi_1\otimes \pi_2}^1\otimes \id_{\rho_1}
 \otimes \id_{\rho_2}(w(\mathbf{1}\hat{\otimes }\rho_1, \mathbf{1}\hat{\otimes}\rho_2)).
 \end{eqnarray*}

 (We identify $H_{\pi_1\hat{\otimes }\rho_1}\otimes H_{\pi_2\hat{\otimes }\rho_2}$
 and $H_{\pi_1}\otimes H_{\pi_2}\otimes H_{\rho_1}\otimes H_{\rho_2}
 $ in the canonical way.)

 We first verify $ w(\pi\hat{\otimes}\mathbf{1},\rho\hat{\otimes}\mathbf{1}) \in
 M^\omega\otimes B(H_\pi)\otimes B(H_\rho)$.
 This follows from the following computation. Here we extend $V_\pi$ and
 $U_\pi$ for $\pi\in \rG$ as in the remark in \S \ref{sub:inter}. 
 \begin{eqnarray*}
  w(\pi\hat{\otimes}\mathbf{1},\rho\hat{\otimes}\mathbf{1}) 
 &=&
  V_\pi^{12}
  V_{\rho}^{13}
 V_{\pi\otimes \sigma}^*\\  
 &=&
 \sum_{\sigma,a} U_\pi^{12*}\lambda_\pi^{12}
  U_{\rho}^{13*}\lambda_{\rho}^{13}
 T_{\pi,\rho}^{\sigma,a} \lambda_\sigma^*U_\sigma T_{\pi,\rho}^{\sigma,a*} \\  
 &=&
 \sum_{\sigma,a} U_\pi^{12*}
  \alpha_\pi^\omega\otimes \id_\rho(U_{\rho}^*)\lambda_\pi^{12}
 \lambda_{\rho}^{13}
 T_{\pi,\rho}^{\sigma,a} \lambda_\sigma^*U_\sigma T_{\pi,\rho}^{\sigma,a*} \\  
 &=&
 \sum_{\sigma,a} U_\pi^{12*}
  \alpha_\pi^\omega\otimes \id_\rho(U_{\rho}^*)
 T_{\pi,\rho}^{\sigma,a} U_\sigma T_{\pi,\rho}^{\sigma,a*} \\
 &=&U_{\pi}^{12*}
  \alpha_\pi^\omega\otimes \id_\rho(U_{\rho}^*)
 U_{\pi\otimes \rho} \\
 &\in & M^\omega\otimes B(H_\pi)\otimes B(H_\rho).  
 \end{eqnarray*}

 Let us examine $
 w(\pi\hat{\otimes}\mathbf{1},\rho\hat{\otimes}\mathbf{1}) $ commutes with 
 $x\otimes 1_\pi\otimes 1_\rho \in M\otimes \mathbf{C}1_\pi\otimes
 \mathbf{C}1_\rho$.
 Note that $\Ad
 U_\pi(x)=\alpha_\pi(x)$ holds for $\pi\in \rG$, $x\in M$. Thus
 \begin{eqnarray*}
 w(\pi\hat{\otimes}\mathbf{1},\pi\hat{\otimes}\mathbf{1})
 (x\otimes 1_\pi\otimes 1_\rho)&=&
  U_\pi^{12*}
  \alpha_\pi^\omega\otimes \id_\rho(U_{\rho}^*)
  U_{\pi\otimes \rho} (x\otimes 1_\pi\otimes 1_\rho) \\
 &=& 
  U_\pi^{12*}
  \alpha_\pi^\omega\otimes \id_\rho(U_{\rho}^*)(\alpha_\pi\otimes \id_\rho)\circ \alpha_\rho(x)
 U_{\pi\otimes \rho} \\
 &=& 
  U_\pi^{12*}
  \alpha_\pi^\omega\otimes \id_\rho(U_{\rho}^*\alpha_\rho(x))
 U_{\pi\otimes \rho} \\
 &=& 
  U_\pi^{12*}
  \alpha_\pi^\omega\otimes \id_\rho(x\otimes 1_\rho)\alpha_\pi\otimes \id_\rho(U_\rho^*)
 U_{\pi\otimes \rho} \\
 &=& (x\otimes 1_\pi\otimes 1_\rho)
 U_\pi^{12*}
  \alpha\otimes \id_\rho(U_\rho^*)
 U_{\pi\otimes \rho} \\
 &=&
 (x\otimes 1_\pi\otimes 1_\rho)
 w(\pi\hat{\otimes}\mathbf{1},\rho\hat{\otimes}\mathbf{1})
 \end{eqnarray*}
 holds, and $w(\pi\hat{\otimes }\mathbf{1},\rho\hat{\otimes}\mathbf{1})\in
 M_\omega\otimes B(H_\pi)\otimes B(H_\rho)$.

 It is trivial that $w(\mathbf{1}\hat{\otimes }\pi,
 \mathbf{1}\hat{\otimes}\rho)=U_\pi^{12*}U_\rho^{13*}U_{\pi\otimes \rho} 
 \in M^\omega \otimes B(H_\pi)\otimes
 B(H_\rho)$. Thus we only have to see 
 $[w(\mathbf{1}\hat{\otimes }\pi, \mathbf{1}\hat{\otimes}\rho), (x\otimes
 1_\pi\otimes 1_\rho)]=0$, $x\in M$. This follows from the following
 computation. 
 \begin{eqnarray*}
 w(\mathbf{1}\hat{\otimes }\pi, \mathbf{1}\hat{\otimes}\rho)
 (x\otimes 1_\pi\otimes 1_\rho)&=&
 U_\pi^{12*}U_{\rho}^{13*}
 U_{\pi\otimes \rho} (x\otimes 1_\pi\otimes 1_\rho) \\
 &=&U_\pi^{12*}U_{\rho}^{13*}
 (\alpha_\pi\otimes \id_\rho)\circ \alpha_\rho(x)
 U_{\pi\otimes \rho} \\
 &=&
 U_\pi^{12*}U_{\rho}^{13*}
 \left(F_{\rho,\pi}\alpha_\rho\otimes \id_\pi(\alpha_\pi(x))F_{\pi,\rho}\right)
 U_{\pi\otimes\rho} \\
 &=&(x\otimes 1_\pi\otimes 1_\rho)
 U_\pi^{12*}U_{\rho}^{13*}
 U_{\pi\otimes\rho} \\
 &=&(x\otimes 1_\pi\otimes 1_\rho)
 w(\mathbf{1}\hat{\otimes }\pi, \mathbf{1}\hat{\otimes}\rho),
 \end{eqnarray*}
 where we used the commutativity of $\alpha_\pi$ and $\alpha_\rho$ in the
 third equality.

 Finally we verify $ w(\mathbf{1}\hat{\otimes }\pi, \rho\hat{\otimes
 }\mathbf{1})\in M_\omega\otimes B(H_\rho)\otimes B(H_\pi)$. As in the
 above, we will see 
 $ w(\mathbf{1}\hat{\otimes }\pi, \rho\hat{\otimes }\mathbf{1})\in
 M^\omega \otimes B(H_\rho)\otimes B(H_\pi)$ and 
 $[ w(\mathbf{1}\hat{\otimes }\pi, \rho\hat{\otimes }\mathbf{1}),
 x\otimes 1_\rho\otimes 1_\pi]=0$, $x\in M$ separately as follow. 
 \begin{eqnarray*}
  w(\mathbf{1}\hat{\otimes }\pi, \rho\hat{\otimes }\mathbf{1})&=& 
  U_\pi^{13*}V_\rho^{12}U_\pi^{13}V_\rho^{12*}\\
 &=&
 U_\pi^{13*}U_\rho^{12*}\lambda_\rho^{12}U_\pi^{13}
 \lambda_\rho^{12*}U_\rho^{13}\\
 &=&
 U_\pi^{13*}U_\rho^{12*}\alpha_\rho^\omega\otimes \id_\pi(U_\pi)
 U_\rho^{12}\\
 &\in &M^\omega\otimes B(H_\rho)\otimes B(H_\pi).
 \end{eqnarray*}

 \begin{eqnarray*}
  w(\mathbf{1}\hat{\otimes }\pi, \rho\hat{\otimes }\mathbf{1})(x\otimes 1_\rho\otimes 1_\pi)
 &=& 
 U_\pi^{13*}U_\rho^{12*}\alpha_\rho^\omega\otimes\id_\pi(U_\pi)
 U_\rho^{12}(x\otimes 1_\rho\otimes 1_\pi) \\
 &=& 
 U_\pi^{13*}U_\rho^{12*}(\alpha_\rho^\omega\otimes \id_\pi)
 (U_\pi)(\alpha_\rho(x)\otimes 1_\pi)
 U_\rho^{12}\\
 &=& 
 U_\pi^{13*}U_\rho^{12*}(\alpha_\rho^\omega\otimes 1_\pi)(U_\pi
 (x\otimes 1_\pi))
 U_\rho^{12}\\
 &=& 
 U_\pi^{13*}U_\rho^{12*}(\alpha_\rho^\omega\otimes 1_\pi\circ \alpha_\pi(x))
 (\alpha_\rho\otimes \id_\pi(U_\pi))
 U_\rho^{12}\\
 &=& (x\otimes 1_\pi\otimes 1_\rho)
 U_\pi^{13*}U_\rho^{12*}
 \alpha_\rho^\omega\otimes \id_\pi(U_\pi)
 U_\rho^{12}\\
 &=& (x\otimes 1_\rho\otimes 1_\pi)
  w(\mathbf{1}\hat{\otimes }\pi, \rho\hat{\otimes }\mathbf{1}).
 \end{eqnarray*}

 \hfill$\Box$

 \begin{lem}\label{lem:repomega}
 We can choose $\bar{U}_\pi\in M^\omega$ 
 so that $\Ad \bar{U}_\pi(x\otimes 1)=\alpha_\pi(x)$, $x\in M$, and 
 $\bar{U}_\pi$ and $\bar{U}_\pi^*\lambda_\pi$ are both
 representation of $\hG$ with $[\bar{U}_{\pi_{ij}}, (\bar{U}_\rho^*\lambda_\rho)_{kl}]=0$.
 \end{lem}
 \textbf{Proof.} 
 By Lemma \ref{lem:cocycleaction}, $\gamma_{\pi\hat{\otimes}\rho}$ defines a cocycle
 twisted action of $\hG\times \hG$ on $M_\omega$. Hence by Theorem
 \ref{thm:2vanish-1}, $w(\cdot,\cdot)=\partial_{\gamma}(u)(\cdot,\cdot)$
 for some $u_{\pi\hat{\otimes}\rho}\in M_\omega\otimes B(H_{\pi\hat{\otimes}\rho})$.  
 Set $\tilde{U}_\pi=U_\pi u_{\mathbf{1}\hat{\otimes}\pi}^*$, and
 $\tilde{V}_\pi=u_{\pi\hat{\otimes}\mathbf{1}}
 V_\pi$.
 By the remark after Definition \ref{df:2coboundary}, 
 $u_{\pi\hat{\otimes}\rho}V^{12}_{\pi}U_\rho^{13*}$ is a
 1-cocycle for the trivial action of $\hG\otimes \hG$.
 This implies that
 $\tilde{U}_\pi^*$ and $\tilde{V}_\pi$ are 
 representation of $\hG$ in $M_1^\omega$ with
 $[\tilde{U}_{\pi_{ij}}, 
 \tilde{V}_{\rho_{kl}}]=0$.
 Moreover we have 
 $$\Ad \tilde{U}_\pi(x)=\Ad U_{\pi}u_{\mathbf{1}\hat{\otimes}\pi}^*(x\otimes 1)=\Ad
 U_\pi(x\otimes 1)=\alpha_\pi(x), x\in M.$$

 Put
 $c_\pi:=\tilde{V}_\pi\lambda_\pi^*\tilde{U}_\pi=u_{\pi\hat{\otimes}\mathbf{1}}
 U_\pi^*\lambda_\pi
 \lambda_\pi^*U_\pi
 u_{\mathbf{1}\hat{\otimes}\pi}^*=u_{\pi\hat{\otimes}\mathbf{1}}
 u_{\mathbf{1}\hat{\otimes}\pi}^*\in M_\omega$.

 Define $W_{\pi_{ij}}:=\sum_k\tilde{V}_{\pi_{ik}}\tilde{U}_{\pi_{kj}}$,
 Set
 $P:=M_\omega\vee\{W_{\pi_{ij}}\}(\subset M_1^\omega)$. 
 Then $M_\omega \subset P$ is the quantum double
 for $\Ad u_\pi\gamma_\pi$. 
 (Unitaries $\tilde{V}_\pi$, $\tilde{U}_\pi$, and $W_\pi$ correspond to $v_\pi$, $u_\pi$,
 and $w_\pi$ in the proof of Lemma \ref{lem:LRunitary}.)
 We have the (unique) conditional expectation
 $E$ from $P$ on $M_\omega$, and it satisfies $E(W_{\pi_{ij}})=\delta_{\pi,\mathbf{1}}$.

 We next prove $\sum_i\lambda_{\pi_{ii}}\in P$.
 Since $c_\pi=\tV_\pi \lambda_\pi^*\tU_\pi$, we
 have $\lambda_{\pi_{ij}}=\sum_{k,l}\tU_{\pi_{ik}}c_{lk}^*\tV_{\pi_{lj}}$.
 Note $\gamma_\pi^2(x):=\tU_\pi^*(x\otimes 1)\tU_\pi$ is an action on
 $M_\omega$, and it follows
 $\tU_{\pi_{ji}}^*x=\sum_{k}\gamma_\pi^2(x)_{ik}\tU_{\pi_{jk}}^*$, $x\in M_\omega$.

 Then we have
 \begin{eqnarray*}
  \sum\lambda_{\pi_{ii}}&=&
 \sum_{i,k,l}\tU_{\pi_{ik}}c_{lk}^*\tV_{\pi_{li}} \\
 &=&
 \sum_{i,k,l}\tU_{\bar{\pi}_{ik}}^*c_{lk}^*\tV_{\pi_{li}} \\
 &=&\sum_{i,j,k,l} \gamma^2_{\bar{\pi}}(c_{lk}^*)_{kj}\tU_{\bar{\pi}_{ij}}^*\tV_{\pi_{li}} \\
 &=&\sum_{i,j,k,l} \gamma^2_{\bar{\pi}}(c_{lk}^*)_{kj}\tU_{\pi_{ij}}\tV_{\pi_{li}} \\
 &=&\sum_{j,k,l} \gamma^2_{\bar{\pi}}(c_{lk}^*)_{kj}W_{\pi_{lj}},
 \end{eqnarray*}
 and thus $\sum_{i}\lambda_{\pi_{ii}}\in P$.

 Define $e:=|G|^{-1}\sum_{\pi,i,j}d\pi W_{\pi_{ii}}$, and 
 $f:=|G|^{-1}\sum_{\pi,i,j}d\pi \lambda_{\pi_{ii}}$. Then $e$ and $f$ are
 projections in $P$
 with $E(e)=E(f)=1/|G|$. 
 Hence there exists a unitary $z\in M_\omega$ such that
 $zfz^*=e$ by Lemma \ref{lem:push}.
 (Though $P$ is not a crossed product of $M_\omega$ by $\hG$, the proof
 of Lemma \ref{lem:push} works for $M_\omega\subset P$ since
 $\{W_\pi\}$ is a representation of $\hG$, and $a\in P$ can be expressed as 
 $\sum_{\pi,i,j}a_{\pi,i,j}W_{\pi_{ij}}$.)

 On one hand, we have 
 \begin{eqnarray*}
  |G|zfz^*&=&\sum_{\pi,i}d\pi z\lambda_{\pi_{ii}}z^* \\
 &=& \sum_{\pi,i,j}d\pi z\alpha_{\pi}(z^*)_{ij}\lambda_{\pi_{ji}}. 
 \end{eqnarray*} 

 On the other hand, since $\tV_\pi=c_\pi \tU_\pi^*\lambda_\pi,$we get
 \begin{eqnarray*}
  |G|e&=&\sum_{\pi,i,j}d\pi \tU_{\pi_{ij}}\tV_{\pi_{ji}} \\
  &=&\sum_{\pi,i,j,k,l}d\pi \tU_{\pi_{ij}}c_{\pi_{jk}}\tU_{\pi_{lk}}^*\lambda_{li}.
 \end{eqnarray*}
 Since $z\alpha_\pi(z^*)_{ij}, \tU_{\pi_{ij}}c_{\pi_{jk}}\tU_{\pi_{lk}}^*\in
 M^\omega$, we have 
 $z\alpha_\pi(z^*)_{il}= \sum_{j,k}\tU_{\pi_{ij}}c_{\pi_{jk}}\tU_{\pi_{lk}}^*$, and
 this implies $(z\otimes 1)\alpha_{\pi}^\omega(z^*)=\tU_\pi c_\pi \tU_\pi^*$. 

 Define $\bar{V}_\pi$ and $\bar{U}_\pi$ by 
 $\bar{V}_\pi=(z^*\otimes 1)\tV_\pi(z\otimes 1)$, 
 $\bar{U}_\pi=(z^*\otimes 1)\tU_\pi(z\otimes 1)$.
 We have 
 \begin{eqnarray*}
  \bar{V}_\pi\lambda_\pi^*\bar{U}_\pi&=&(z^*\otimes 1)\tV_\pi(z\otimes 1)\lambda_\pi^*
 (z^*\otimes 1)\tU_\pi(z\otimes 1) \\
 &=&(z^*\otimes 1)\tV_\pi\lambda_\pi\alpha_\pi(z)
 (z^*\otimes 1)\tU_\pi(z\otimes 1) \\
 &=&(z^*\otimes 1)\tV_\pi\lambda_\pi \tU_\pi c_\pi^* \tU_\pi^*
 \tU_\pi(z\otimes 1) \\
 &=&1,
 \end{eqnarray*}
 hence $\bar{V}_\pi=\bar{U}_\pi^*\lambda_\pi$. 
 Since $z\in M_\omega$, $\Ad \bar{U}_\pi(x\otimes 1)=\alpha_\pi(x)$
 holds for $x\in
 M\subset M^\omega$. It is clear that $\bar{U}_\pi$ and $\bar{V}_\pi$ are both
 representations of $\hG$ with $[\bar{V}_{\pi_{ij}},\bar{U}_{\rho_{kl}}]=0$.  
 \hfill$\Box$

 \noindent
 \textbf{Remark.} 
 To avoid using the commutativity of $\hG$,
 we should consider $\gamma_\pi^3(x):=\Ad U_{\bar{\pi}}^*(x\otimes 1)$
 rather than $\gamma_\pi^2=\Ad U_{\pi}^*(x\otimes 1)$. By suitable inner
 perturbation, this $\gamma_\pi^3$ is shown to be an
 ``conjugate'' action of $\hG$ in the sense $(\gamma_\pi^3\otimes
 \id_\rho)\circ\gamma^3_\rho(x)\overline{T}=\overline{T}\gamma_\sigma^3(x) $ for
 $T\in (\sigma,\pi\otimes \rho)$ without using the commutativity of $\hG$.  
 (See the remark after Lemma \ref{lem:reprel}.)

 \begin{cor}\label{cor:repomega}
 Fix $U_\pi$ as in Lemma \ref{lem:repomega}. Then we have
 $\alpha_\rho^\omega(U_{\pi_{ij}})=\Ad U_\rho(U_{\pi_{ij}}\otimes 1_\rho)$. 
 \end{cor}
 \textbf{Proof.} By Lemma \ref{lem:repomega}, we have $U_\pi V_\pi=\lambda_\pi$. Since
 $[U_{\pi_{ij}}, V_{\rho_{kl}}]=0$, we have
 $$\alpha^\omega_\rho(U_{\pi_{ij}})=
 \Ad \lambda_\pi(U_{\pi_{ij}}\otimes1_\rho)=
 \Ad U_\rho V_\rho(U_{\pi_{ij}}\otimes1_\rho)=
 \Ad U_\rho (U_{\pi_{ij}}\otimes1_\rho).$$
 \hfill$\Box$

 \begin{lem}\label{lem:muomega}
 We choose $U_\pi$ as in Lemma \ref{lem:repomega}.
  There exists an $\alpha^\omega$-equivariant system of matrix units
  $E=\{E_{\pi_{ij},\rho_{kl}}\}\subset M^\omega$ such that
  $\lambda_{\pi}^E=U_\pi$ and $E_{\mathbf{1},\mathbf{1}}\in M_\omega$.
 \end{lem}
 \textbf{Proof.} 
 Let $\gamma^i_\pi$ and
 $\gamma_{\pi\hat{\otimes}\rho}=(\gamma_\pi^1\otimes
 1_\rho)\circ\gamma^2_\rho(x)$ be as in Lemma \ref{lem:cocycleaction}. By Lemma
 \ref{lem:repomega}, $\gamma$ is an outer action of $\hG\times \hG$ on $M_\omega$.
 By Corollary \ref{cor:equmu}, there
 exists a $\gamma$-equivariant system of matrix units
 $\{e_{(\xi_{ab}\hat{\otimes}\pi_{ij}), (\eta_{cd}\hat{\otimes}\rho_{kl})}\}$
 in $M_\omega$. Put
 $F_{\pi_{ij},\rho_{kl}}:=\sum\limits_{\xi,a,b}
 e_{(\xi_{ab}\hat{\otimes}\pi_{ij}), (\xi_{ab}\hat{\otimes}\rho_{kl})}$.
 Then $\{F_{\pi_{ij},\rho_{kl}}\}$ is in $M_\omega^{\gamma^1}$, and
 becomes a 
 $\gamma^2$-equivariant system of matrix units.

 Set
 $\tilde{F}^\pi_{i,j}:=\sum_{k}F_{\pi_{ik},\pi_{jk}}$.
 Then it is easy to see that $\{\tilde{F}_{i,j}^\pi\}$ forms a system of
 matrix units for $R(G)$. Namely we have
 $\tilde{F}_{i,j}^{\pi}\tilde{F}_{k,l}^\rho=\delta_{\pi,\rho}\delta_{j,k}\tilde{F}_{i,l}^\pi$, 
 $\tilde{F}_{i,j}^{\pi*}=\tilde{F}_{j,i}^\pi$ and $\sum_{\pi,i}\tilde{F}_{i,i}^\pi=1$.
 Since $F=\{F_{\pi_{ij},\rho_{kl}}\}$ is $\gamma^2$-equivariant,
 we have $\gamma_{\pi}^2(\tilde{F}^{1})_{i,j}=\sum_k
  \lambda^F_{\pi_{ik}}F_{\mathbf{1},\mathbf{1}}
 \lambda^F_{\bar{\pi}_{jk}}=d\pi^{-1} \tilde{F}_{i,j}^\pi$.
 Since $\gamma^2_\pi(x)=\Ad U_\pi^*(x\otimes 1_\pi)$ by definition and 
 $U_{\pi_{ij}}^*=U_{\bar{\pi}_{ij}}$, we
 have $U_{\bar{\pi}_{ji}}x=\sum_k\gamma_\pi^2(x)_{ik}U_{\bar{\pi}_{jk}}$.

 Define $E_{\pi_{ij},\rho_{kl}}:=\sqrt{d\pi
 d\rho}U_{\pi_{ij}}F_{\mathbf{1},\mathbf{1}}
 U_{\bar{\rho}_{kl}}$.
 It is trivial $E_{\mathbf{1},\mathbf{1}}\in M_\omega$.

 We first prove $\lambda_{\pi}^E=U_{\pi}$.
 \begin{eqnarray*}
  \lambda^E_{\pi_{ij}}&=&\sum_{\rho,{k.l},\sigma,{m,n},e}\sqrt{\frac{d\rho}{d\sigma}}
 T_{\pi_i,\rho_k}^{\sigma_m,e} \overline{T_{\pi_j,\rho_l}^{\sigma_n,e}}
 E_{\sigma_{mn},\rho_{kl}} \\
 &=&
 \sum_{\rho,{k,l},\sigma,{m,n},e}{d\rho}
 T_{\pi_i,\rho_k}^{\sigma_m,e}U_{\sigma_{mn}} \overline{T_{\pi_j,\rho_l}^{\sigma_n,e}}
 F_{\mathbf{1},\mathbf{1}}
 U_{\bar{\rho}_{kl}} \\
 &=&
 \sum_{\rho,{k,l}}d\rho
 U_{\pi_{ij}} U_{\rho_{kl}}
 F_{\mathbf{1},\mathbf{1}}
 U_{\bar{\rho}_{kl}} \\
 &=&
 \sum_{\rho,k,l,m}d\rho
 U_{\pi_{ij}} \gamma_{\bar{\rho}}^2(F_{\mathbf{1},\mathbf{1}})_{lm}U_{\rho_{km}}
 U_{\bar{\rho}_{kl}}\\
 &=&
 \sum_{\rho,l,m}d\rho
 U_{\pi_{ij}} \gamma_{\bar{\rho}}^2(F_{\mathbf{1},\mathbf{1}})_{lm}
 \left(\sum_k U_{\bar{\rho}_{km}}^*
 U_{\bar{\rho}_{kl}}\right)\\
 &=&
 \sum_{\rho,l}d\rho
 U_{\pi_{ij}} \gamma_{\bar{\rho}}^2(F_{\mathbf{1},\mathbf{1}})_{ll}\\
 &=&\sum_{\rho,l}U_{\pi_{ij}}\tilde{F}_{l,l}^\rho \\
 &=&U_{\pi_{ij}}.
 \end{eqnarray*}

 We next prove $E=\{E_{\pi_{ij},\rho_{kl}}\}$ is a system of matrix units.

 If we set $\pi =\mathbf{1}$ in the above computation, 
 we get $\sum_{\pi,i,j} E_{\pi_{ij},\pi_{ij}}=1$. 
 It is easy to see $E_{\pi_{ij},\rho_{kl}}^*
 =E_{\rho_{kl},\pi_{ij}}$. Thus we only have to verify 
 $E_{\pi_{ij},\rho_{kl}}E_{\sigma_{mn},\xi_{ab}}=
 \delta_{\rho,\sigma}\delta_{k,m}\delta_{l,n}
 E_{\pi_{ij},\xi_{ab}}$. At first we compute 
 $ F_{\mathbf{1},\mathbf{1}}U_{{\pi}_{ij}}U_{\rho_{kl}}F_{\mathbf{1},\mathbf{1}}$.
 Note 
 $F_{\mathbf{1},\mathbf{1}}\gamma_\sigma(F_{\mathbf{1},\mathbf{1}})_{m,n}=
 d\sigma^{-1}\tilde{F}^\mathbf{1}\tilde{F}^\sigma_{m,n}=
 \delta_{\mathbf{1},\sigma}F_{\mathbf{1},\mathbf{1}}$. 

 Then 
 \begin{eqnarray*}
  F_{\mathbf{1},\mathbf{1}}U_{{\pi}_{ij}}U_{\rho_{kl}}F_{\mathbf{1},\mathbf{1}}
 &=&\sum_{\sigma,m,n}
  F_{\mathbf{1},\mathbf{1}}T_{{\pi}_{i},\rho_k}^{\sigma_m}U_{\sigma_{mn}}
 \overline{T_{\pi_j,\rho_l}^{\sigma_n}}F_{\mathbf{1},\mathbf{1}} \\
 &=&\sum_{\sigma,m,n,a}
 T_{{\pi}_{i},\rho_k}^{\sigma_m}\overline{T_{\pi_j,\rho_l}^{\sigma_n}}
 F_{\mathbf{1},\mathbf{1}}
 \gamma^2_{\bar{\sigma}}(F_{\mathbf{1},\mathbf{1}})_{m,a} U_{\sigma_{na}} \\
 &=&
 \sum_{\sigma,m,n,a}d\sigma^{-1}
 T_{{\pi}_{i},\rho_k}^{\sigma_m}\overline{T_{\pi_j,\rho_l}^{\sigma_n}}
 \tilde{F}^{\mathbf{1}}
 \tilde{F}^{\bar{\sigma}}_{m,a} U_{\sigma_{na}} 
 \\
 &=&
 T_{{\pi}_{i},\rho_k}^{\mathbf{1}}\overline{T_{\pi_j,\rho_l}^{\mathbf{1}}}
 \tilde{F}^{\mathbf{1}}
 \\
 &=&\delta_{\pi,\bar{\rho}}\delta_{i,k}\delta_{j,l}d\pi^{-1}F_{\mathbf{1},\mathbf{1}}
 \end{eqnarray*}
 holds, and hence we have
 \begin{eqnarray*}
 E_{\pi_{ij},\rho_{kl}}E_{\sigma_{mn},\xi_{ab}}
 &=&
 \sqrt{d\pi d\rho d\sigma d\xi} 
 U_{\pi_{ij}}F_{\mathbf{1},\mathbf{1}}
 U_{\bar{\rho}_{kl}}U_{\sigma_{mn}}F_{\mathbf{1},\mathbf{1}}U_{\bar{\xi}_{ab}} \\
  &=&
 \delta_{\rho,\sigma}\delta_{k,m}\delta_{l,n}d\rho^{-1}
 \sqrt{d\pi d\rho d\sigma d\xi} 
 U_{\pi_{ij}}F_{\mathbf{1},\mathbf{1}}
 U_{\bar{\xi}_{ab}} \\
 &=&
 \delta_{\rho,\sigma}\delta_{k,m}\delta_{l,n}
 E_{\pi_{ij},\xi_{ab}}.
 \end{eqnarray*}

 Finally we verify that $\{E_{\pi_{ij},\rho_{kl}}\}$ is an
 $\alpha^\omega$-equivariant system of matrix units.
 Since $F_{\mathbf{1},\mathbf{1}}\in M_\omega^{\gamma^1}$, 
 $\alpha_\pi^\omega(F_{\mathbf{1},\mathbf{1}})=\Ad
 U_\pi(F_{\mathbf{1},\mathbf{1}}\otimes 1)$ holds.
 Together with
  Corollary \ref{cor:repomega},
 we have 
 \begin{eqnarray*}
  \alpha_\sigma^\omega(E_{\pi_{ij},\rho_{kl}})&=&\sqrt{d\pi d\rho}\,
 \alpha_\sigma^\omega(U_{\pi_{ij}})\alpha_{\sigma}^\omega(F_{\mathbf{1},\mathbf{1}})
 \alpha_\sigma^\omega(U_{\bar{\rho}_{kl}}) \\
 &=&\sqrt{d\pi d\rho}\,
 \Ad U_\sigma(U_{\pi_{ij}}\otimes 1_\sigma)(F_{\mathbf{1},\mathbf{1}}
 \otimes 1_\sigma)(U_{\bar{\rho}_{kl}}
 \otimes 1_\sigma) \\
 &=&\Ad U_\sigma(E_{\pi_{ij},\rho_{kl}}\otimes 1_\sigma) \\
 &=&\Ad \lambda_\sigma^E(E_{\pi_{ij},\rho_{kl}}\otimes 1_\sigma).
 \end{eqnarray*}
 \hfill$\Box$

 \noindent
 \textbf{Remark.} We can regard $\{\tilde{F}^\pi_{i,j}\}$ as an
 analogue Rohlin projections for $\gamma^2$.

 \begin{prop}\label{prop:lift}
  Let $E=\{E_{\pi_{ij},\rho_{kl}}\}\subset M^\omega$ be an
  $\alpha^\omega$-equivariant system of 
  matrix units. Then there exists a representing sequence of systems of matrix units
  $\{e_{\pi_{ij},\rho_{kl}}^n\}$ for $E_{\pi_{ij},\rho_{kl}}$, and
  1-cocycles $\{u_\pi^n\}$ for $\alpha$, $n=1,2,3\cdots$, such that
  $(u_\pi^n)=1$ in $M^\omega$ and 
 each $\{e_{\pi_{ij},\rho_{kl}}^n\}$ is $\Ad u_\pi^n\alpha_\pi$-equivariant.
  \end{prop}
 \textbf{Proof.} 
 Fix a representing sequence $\{e_{\pi_{ij},\rho_{kl}}^n\}$ for
 $E_{\pi_{ij},\rho_{kl}}$ consisting of systems of matrix units. Set
 $A_n:=\{e_{\pi_{ij},\rho_{kl}}^n\}''\subset M $, 
 and $\lambda^n_\pi$ the unitary representation of $\hG$ associated with $A_n$.
 Since $E_{\mathbf{1},\mathbf{1}}\otimes 1=
 (E_{\mathbf{1},\mathbf{1}}\otimes
 1)\lambda_\pi^{E*}\lambda_\pi^E(E_{\mathbf{1},\mathbf{1}}
 \otimes 1)$ and $\alpha_\pi(E_{\mathbf{1},\mathbf{1}})=\Ad
 \lambda_\pi^E( E_{\mathbf{1},\mathbf{1}}\otimes 1)=
 \lambda_\pi^E(E_{\mathbf{1},\mathbf{1}}\otimes
 1)(E_{\mathbf{1},\mathbf{1}}\otimes 1)
 \lambda_\pi^{E*}$, we can choose a representing sequence $\{v_\pi^n\}$ for 
 $(E_{\mathbf{1},\mathbf{1}}\otimes 1)\lambda_\pi^*$ such that
 $v_{\pi}^nv_\pi^{n*}=e_{\mathbf{1},\mathbf{1}}^n\otimes 1$,
 and $v_{\pi}^{n*}v_\pi^n=\alpha_\pi(e_{\mathbf{1},\mathbf{1}}^n)$ by Lemma \ref{lem:ultra}.
 Set $w_\pi^n:=\sum_{\pi,i,j}(e_{\pi_{ij},\mathbf{1}}^n\otimes
 1)v_\pi^n\alpha_\pi(e^n_{\mathbf{1},\pi_{ij}})$. 
 Then $w_\pi^n$ is a unitary, and 
 $\Ad w_\pi^n\alpha_\pi(e_{\sigma_{ij},\rho_{kl}})
 =(e_{\sigma_{ij},\rho_{kl}}\otimes 1)$
 holds. Define $\alpha_\pi^n:=\Ad w_\pi^n\alpha_\pi$, and 
 $U_{\pi,\rho}^n:=(\partial_{\alpha^n}w^n)_{\pi,\rho}
 $. Since $\alpha_\pi^n$ is trivial on $A_n$,
 $U_{\pi,\rho}^n\in (A_n'\cap M)\otimes
 B(H_\pi)\otimes B(H_\rho)$ and
 $\{\alpha_\pi^n,U_{\pi,\rho}^n\}$ is a cocycle twisted action on
 $A_n'\cap M$.

 We have $(w_\pi^n )=
 \sum_{\pi,i,j}(E_{\pi_{ij},\mathbf{1}}\otimes 1) 
 (E_{\mathbf{1},\mathbf{1}}\otimes 1)\lambda_\pi^{E*}
 \alpha_{\pi}(E_{\mathbf{1},\pi_{ij}})=\lambda_\pi^{E*}$, hence
 $U_{\pi,\rho}^n\rightarrow 1$ as $n\rightarrow\omega$. By Theorem
 \ref{thm:2vanish-2}, there exists $\bar{w}_\pi^n\in U((A_n'\cap M)\otimes B(H_\pi))$ with
 $U_{\pi,\rho}^n=(\partial_{\alpha^n}\bar{w}^{n*})_{\pi,\rho}$
 and $\lim_{n\rightarrow \omega}\|\bar{w}_\pi^n-1\|_2=0$. Set
 $u_\pi^n:=\lambda_\pi^n\bar{w}_\pi^nw_\pi^n$. 
 Then $(u_\pi^n)=(\lambda_\pi^nw_\pi^n)=\lambda_\pi^E\lambda_\pi^{E*}=1$ in $M^\omega$, and 
  $u_\pi^n$ is a 1-cocycle for $\alpha_\pi$ 
 by Lemma \ref{lem:repcocycle} and the remark
  after Definition \ref{df:2coboundary}.
 It is trivial that $\Ad u_\pi^n\alpha_\pi=\Ad
  \lambda_\pi^n$ on $A_n$, and hence $\{e^n_{\pi_{ij},\rho_{kl}}\}$ is
  $\Ad u_\pi^n\alpha_\pi$-equivariant.
 \hfill$\Box$

 \section{Classification}\label{sec:class}
 \begin{prop}\label{prop:idsplit}
  Let $\alpha$ be an outer action on $\cR$. Then $\alpha$ is conjugate to
  $\alpha\otimes \id_{\cR}$.
 \end{prop}
 \textbf{Proof.} This follows from \cite{Bi-Mc} since $\cR'\cap (\cR^\alpha)^\omega$ is
 noncommutative. \hfill$\Box$

 \begin{lem}\label{lem:relcommexp}
  Let $K\subset \cR$ be a subfactor with $K\cong M_n(\mathbf{C})$, and
  $\{e_{ij}\}$ be a system of  matrix units for $K$. If
  $\|[x,e_{ij}]\|_2<\varepsilon/n$,  then $\|E_{K'\cap \cR}(x)-x\|_2<\varepsilon$.
 \end{lem}
 \textbf{Proof.} Since $E_{K'\cap
 \cR}(x)=\frac{1}{n}\sum_{i,j}e_{ij}xe_{ji}$, 
 \begin{eqnarray*}
  \|E_{K'\cap \cR}(x)-x\|_2&\leq& \frac{1}{n}\sum_{i,j}\|e_{ij}xe_{ji}-xe_{ij}e_{ji}\|_2 \\
 &\leq &\frac{1}{n}\sum_{i,j}\|[e_{ij},x]\|_2 \\
 &< & \varepsilon
 \end{eqnarray*}
 holds. \hfill$\Box$

 \begin{lem}\label{lem:modelpiece}
 For any $\varepsilon>0$, $a_1,a_2,\cdots, a_n\in \cR$, there exist a
  1-cocycle $u_\pi$ for $\alpha_\pi$, and an
  $\Ad u_\pi\alpha$-equivariant system of matrix units $E=\{e_{\pi_{ij},\rho_{kl}}\}$ such that 
 $\|\alpha_\pi(a_i)-\Ad \lambda_\pi^E(a_i\otimes 1)\|_2< \varepsilon$, 
 $\|u_\pi -1\|_1<\varepsilon$
 and $\|[e_{\mathbf{1},\mathbf{1}},a_i]\|_2<\varepsilon$
 \end{lem}
 \textbf{Proof.} 
 By Lemma \ref{lem:muomega} and 
 Proposition \ref{prop:lift}, we have systems of matrix units
 $E_n=\{e_{\pi_{ij},\rho_{kl}}^n\}$ and 1-cocycles $u_\pi^n$ for $\alpha_\pi$
 such that $\{e_{\pi_{ij},\rho_{kl}}^n\}$ is $\Ad u_\pi^n\alpha_\pi$
 equivariant, $\alpha_\pi(x)=\lim_{n\rightarrow \omega}\Ad
 \lambda_\pi^{E_n} (x\otimes  1)$, $\lim_{n\rightarrow \omega}\|u_\pi^n- 1\|_2=0$ 
 and $\lim_{n\rightarrow \omega}\|[e_{\mathbf{1},\mathbf{1}}^n,x]\|_2=0$
 for any $x\in M$. Put $E:=E_n$ and
 $u_\pi:=u_\pi^n$ for sufficiently
 large $n$. \hfill$\Box$

 Now we can prove the main theorem of this paper.

 \begin{thm}\label{thm:main}
  Let $\alpha$ be an outer action of $\hat{G}$ on $\cR$. Then $\alpha$ is
  conjugate to the model action $m$. 
 \end{thm}
 \textbf{Proof.} 
 We use notations in Section \ref{sec:model}.
 Let $\{a_i\}_{i=1}^\infty$ be a strongly dense countable subset of
 the unit ball of $\cR$.
 We fix a sequence $\{\varepsilon_n\}$ such that 
 $0<9|G|^3\varepsilon_n\leq 2^{-n}.$ Especially we have $\sum_n\varepsilon_n<\infty$.
 We will construct mutually commuting finite dimensional
 subfactors $K_n\cong M_{|G|}(\mathbf{C})$, 
 unitary 1-cocycles $v_\pi^n$ for $\alpha_\pi$, 
 a unitary 1-cocycle $w_\pi^n$ for $\Ad
 \tilde{\lambda}_\pi^{n-1*}v_\pi^{n-1}\alpha_\pi$
 satisfying the
 following conditions inductively.
 \begin{eqnarray*}
 (1.n) && v_\pi^{n}:=\tilde{\lambda}_\pi^{n-1}w_\pi^{n}\tilde{\lambda}_\pi^{n-1*}v_\pi^{n-1
 }, n\geq 2, \\
 (2.n) && \|w_\pi^n-1\|_2<\varepsilon_n, \\
 (3.n) 
 && \|[a_i,e_{\mathbf{1},\mathbf{1}}^n]\|_2<\varepsilon_n, 1\leq i\leq n, \\
 (4.n) &&\Ad v_\pi^n\alpha_\pi=\Ad
 \tilde{\lambda}_\pi^n \mbox{ on } K_1\vee\cdots\vee K_n, \\
 (5.n) &&\|\Ad v_\pi^n\alpha_\pi(a_i)-m_\pi^n(a_i)\|_2<\varepsilon_n, 1\leq i\leq n.
 \end{eqnarray*}

 By Lemma \ref{lem:modelpiece}, we get a unitary cocycle $w_\pi^1$ for
 $\alpha_\pi$, and an $\Ad w_\pi^1\alpha_\pi$-equivariant system of matrix units 
 $\{e_{\pi_{ij},\rho_{kl}}^1\}$ such that
 $\|[a_1,e_{\mathbf{1},\mathbf{1}}^1]\|_2< \varepsilon_1$,
 $\|w_\pi^1-1\|_2<\varepsilon_1$, 
 and 
 $\|\Ad w_\pi^{1}\alpha_\pi(a_1)-\Ad \lambda_\pi^1(a_1\otimes 1)\|_2<\varepsilon_1$. 
 Let $K_1$ be a finite dimensional subfactor generated by $\{e_{\pi_{ij
 },\rho_{kl}}^1\}$, and set $v_\pi^1:=w_\pi^1$. Then we get the conditions
 (2.1), (3.1), (4.1) and (5.1).

 Suppose that we have done up to the $n$-th step. By $(4.n)$, we
 have $\Ad \tilde{\lambda}_\pi^{n*}v_\pi^n\alpha_\pi=\id$ on $K_1\vee
 \cdots \vee K_n$. Hence $\Ad\tilde{\lambda}_\pi^{n*}v_\pi^n\alpha_\pi$ induces
 an action of $\hG$ on $(K_1\vee\cdots K_n)'\cap\cR$.  Decompose $a_i$ as 
 $a_i=\sum_ib_{ik}e_k$, $b_{ik}\in (K_1\vee\cdots K_n)'\cap\cR$, $e_k\in
 K_1\vee\cdots\vee K_n$. By Lemma \ref{lem:modelpiece}, we get a 
 unitary cocycle $w_{\pi}^{n+1}$ for 
 $\Ad \tilde{\lambda}_\pi^{n*}v_\pi^n\alpha_\pi $, a
 $\Ad w_\pi^{n+1}\tilde{\lambda}_\pi^{n*}v_\pi^n\alpha_\pi $
 -equivariant system matrix units
 $K_{n+1}:=\{e_{\pi_{ij},\rho_{kl}}^{n+1}\}\subset (K_1\vee\cdots K_n)'\cap\cR$,
 such that \\
 $(a.n+1)\,\,\|[e_{\mathbf{1},\mathbf{1}}^n, b_{ik}]\|_2<\delta_{n+1},$ \\
 $(b.n+1)\,\,\|w_\pi^{n+1}-1 \|_2<\varepsilon_{n+1},$ \\
 $(c.n+1)\,\,\|\Ad w_\pi^{n+1}\tilde{\lambda}_\pi^{n*}v_\pi^n\alpha_\pi(b_{ik})-\Ad
 \lambda_{\pi}^{n+1}(b_{ik}\otimes 1) \|_2<\delta_{n+1}$,\\
 for sufficiently small $\delta_{n+1}>0$.
 The condition $(b.n+1)$ is nothing but $(2.n+1)$. 
 If we choose sufficiently enough small $\delta_{n+1}$, then we get
 $(3.n+1)$ and 
 $$(c.n+1)'\,\,\|\Ad w_\pi^{n+1}\tilde{\lambda}_\pi^{n*}v_\pi^n\alpha_\pi(a_i)-\Ad
 \lambda_{\pi}^{n+1}(a_i\otimes 1) \|_2<\varepsilon_{n+1}, 1\leq i\leq n+1$$
 from $(a.n+1)$ and $(c.n+1)$ respectively. Set
 $v_\pi^{n+1}:=\tilde{\lambda}_\pi^nw^{n+1}_\pi\tilde{\lambda}_\pi^{n*}v_\pi^n
 $. Then we get $(1.n+1)$ and $(5.n+1)$. 
 Since
  $\{e_{\pi_{ij},\rho_{kl}}^{n+1}\}\subset (K_1\vee\cdots \vee K_n)'\cap\cR$ is 
 $\Ad w_\pi^{n+1}\tilde{\lambda}_\pi^{n*}v_\pi^n\alpha_\pi$-equivariant, 
 we get $(4.n+1)$, and  
 $K_{n+1}$ commutes with $K_i$, $1\leq i\leq n$.
 Thus we complete induction. 

 We will show $\{v_\pi^n\}$ is a Cauchy sequence. 
 \begin{eqnarray*}
  \|v_\pi^{n+1}-v_\pi^n\|_2&=&
  \|v_\pi^{n}\tilde{\lambda}_\pi^nw_{\pi}^{n+1}\tilde{\lambda}_\pi^{n*}-v_\pi^n\|_2\\
 &=&\|\tilde{\lambda}_\pi^nw_{\pi}^{n+1}\tilde{\lambda}_\pi^{n*}-1\|_2\\
 &<&\varepsilon_{n+1}.
 \end{eqnarray*}
 By the choice of $\varepsilon_n$, $\{v_\pi^n\}$ is Cauchy, and
 hence $\lim_{n\rightarrow \infty}v_\pi^n=v_\pi$ exists.

 We will prove $\|[e_{\pi_{ij},\rho_{kl}}^{n+1}, a_i]\|<\varepsilon_n$,
 $1\leq i\leq n$. By $(5.n)$ and $(5.n+1)$, we get 
 $$\|\Ad v_\pi^{n*}\tilde{\lambda}_\pi^{n}(a_i)
 -\Ad v_\pi^{n+1*}\tilde{\lambda}_\pi^{n+1}(a_i)
 \|<2\varepsilon_n, 1\leq i\leq n.$$
 By the definition of $v_\pi^n$, we get $\|a_i\otimes 1-\Ad
 w_\pi^{n+1*}\lambda_\pi^{n+1}(a_i\otimes 1)\|_2< 2\varepsilon_n$.
 Then 
 \begin{eqnarray*}
 \|[a_i\otimes 1,\lambda_\pi^{n+1}]\|_2&=&
 \|[a_i\otimes 1,w_\pi^{n+1}w_\pi^{n+1*}\lambda_\pi^{n+1}]\|_2 \\
 &\leq &\|[a_i\otimes 1,w_\pi^{n+1}]w_\pi^{n+1*}\lambda_\pi^{n+1}\|_2 +
 \|w_\pi^{n+1}[a_i\otimes 1, w_\pi^{n+1*}\lambda_\pi^{n+1}]\|_2 \\
 &\leq &
 \|[a_i\otimes 1,w_\pi^{n+1}-1]\|_2 + 2\varepsilon_n \\
 &<&4\varepsilon_n.
 \end{eqnarray*}
 Hence we get 
 $\|[\lambda_{\pi_{ij}}^{n+1}, a_i]\|_2<4d\pi\varepsilon_n<4|G|\varepsilon_n$ for $1\leq i\leq
 n$. Then we have
 \begin{eqnarray*}
 \lefteqn{\|[a_i,e_{\pi_{ij},\rho_{kl}}^{n+1}]\|_2}\\
 &=&\sqrt{d\pi d\rho}\|[a_i,\lambda_{\pi_{ij}}^{n+1}e_{\mathbf{1},\mathbf{1}}^{n+1}
 \lambda^{n+1}_{{\bar {\rho}}_{kl}}]\|_2 \\ 
 &\leq & |G|(\|[a_i,\lambda^{n+1}_{\pi_{ij}}]e_{\mathbf{1},\mathbf{1}}^{n+1}\|_2
 +\|\lambda_{\pi_{ij}}^{n+1}[a_i,e_{\mathbf{1},\mathbf{1}}^{n+1}]\lambda_{\bar{\rho}_{kl}}^{n+1}\|_2
 +\|\lambda_{\pi_{ij}}^{n+1}e_{\mathbf{1},\mathbf{1}}^{n+1}
 [a_i,\lambda_{\bar{\rho}_{kl}}^{n+1}]\|_2) \\
 &<&9|G|^2 \varepsilon_n \\
 &<&\frac{1}{2^n|G|}.
 \end{eqnarray*}

 This implies $\|E_{K_{n+1}'\cap \cR}(a_i)-\-a_i\|_2< 1/2^{n}$ for $1\leq
 i\leq n$ by Lemma \ref{lem:relcommexp}.
  Set $K:=\bigvee K_n(\cong \cR)$. 
 By \cite[Lemma 2.3.6]{Con-auto},
 $\cR=K\vee K'\cap
 \cR\cong K\otimes K'\cap \cR$. By $(5.n)$,
 $\Ad v_\pi\alpha_\pi=m_\pi \otimes \id_{K'\cap \cR}$. By Proposition \ref{prop:idsplit}
 $m_\pi\cong m_\pi\otimes \id_{\cR}$, and   
 $\Ad v_\pi \alpha_\pi$ is conjugate to $m_\pi\otimes \id_{K\vee K'\cap
 \cR}=m_\pi\otimes \id_{\cR}\cong m_\pi$. By Proposition \ref{prop:1vanish},
 $\alpha$ is conjugate to $m$. \hfill$\Box$

 It is obvious that Theorem \ref{thm:main00}
 follows immediately from Theorem \ref{thm:main}. \\

 \noindent
 \textbf{Remark.} So far we treat  only actions of $\hG$ for a finite group
 $G$. However we can generalize our theory to outer actions of finite
 dimensional Kac algebras. Difference between $\hG$ and general finite
 dimensional Kac algebras is the commutativity $\pi\otimes\rho\cong
 \rho\otimes \pi$. We do not use the commutativity of $\hG$ in proofs
 except Lemma \ref{lem:cocycleaction}. To generalize Lemma
 \ref{lem:cocycleaction} to a finite dimensional Kac algebra
 $\mathcal{K}$, we should consider a (cocycle) 
 action of $\mathcal{K}\otimes \mathcal{K}^\mathrm{opp}$ on $M^\omega$ as in the remark
 after Lemma \ref{lem:repomega}. 
 (Note $R(G)$ and $R(G)^\mathrm{opp}$ are
 essentially same Kac algebras due to cocommutativity of $R(G)$.)

 \appendix
 \section{Twisted crossed product construction}\label{sec:twist}
 Let $\{\alpha, U\}$ be a cocycle twisted action of $\hG$ on $M$. In this
 appendix, we give the definition of a twisted crossed product
 $M\rtimes_{\alpha, U}\hG$.

 Let $H:=L^2(M)$ be the standard Hilbert space. We identify $H\otimes
 \ell^2(\hG)$ with $\{\bigoplus_\pi v(\pi)\mid v(\pi)\in H\otimes B(H_\pi)\}$ as usual.
 Put $\langle v,w \rangle_\pi=\sum_{ij}\langle v_{ij}, w_{ij}\rangle $
 for $v,w \in H\otimes B(H_\pi)$.  
 Then the inner product is given by $\langle v,w \rangle=\sum_{\pi}
 d\pi\langle v(\pi), \omega(\pi)\rangle_\pi $ for $v,w\in H\otimes \ell^2(\hG)$.

 We define an action $\alpha$ of
 $M$ on $H\otimes \ell^2(\hG)$, and
 $\lambda_{\pi_{ij}}\in B(H\otimes \ell^2(\hG))$ 
 by
 $$(\alpha(a)v)(\pi)=\alpha_\pi(a)v(\pi),$$
 $$(\lambda_{\pi_{ij}}v)(\rho):=\sum_{\sigma}U_{\rho,\pi_i}^{\sigma,e}v(\sigma)
 T_{\rho,\pi_j}^{\sigma,e*}.$$

 \begin{df}
  Define $M\rtimes_{\alpha,U}\hG:=\alpha(M)\vee\{\lambda_{\pi_{ij}}\}$,
  and call it the twisted crossed product of $M$ by $\{\alpha, U\}$.
 \end{df}

 \begin{lem} \label{lem:reltwist}
  Set $\lambda_\pi=(\lambda_{\pi_{ij}})\in B(H\otimes \ell^2(\hG))\otimes B(H_\pi)$. 
 Then $\lambda_\pi$ is a unitary, and we have 
  $\lambda_\pi(a\otimes 1_\pi)\lambda_\pi^*=\alpha_\pi(a)$, and
  $\lambda^{12}_\pi\lambda_{\rho}^{13}T_{\pi,\rho}^{\sigma,e}=
 U_{\pi,\rho}^{\sigma,e}\lambda_\sigma$.  
 Set $\tilde{U}_{\pi_i, \bar{\pi}_j}:=\sum_{k}U_{\pi_{ik},\bar{\pi}_{jk}}$. Then 
 we have $\lambda_{\bar{\pi}_{ij}}^*=\sum_kU_{\pi_k,\bar{\pi}_i}^*\lambda_{\pi_{kj}}$.
 Here we identify $\alpha(a)$ and $a$.
 We call $\lambda_\pi$ an implementing unitary.
 \end{lem}

 To show Lemma \ref{lem:reltwist}, we prepare the following lemma.
 \begin{lem}\label{lem:u}
  We have
  $\sum_{k,l}\tilde{U}_{\pi_k,\bar{\pi}_l}^*\alpha_{\pi}
 (\tilde{U}_{\bar{\pi}_l,\pi_i})_{k,j}=\delta_{i,j}$. 
 \end{lem}
 \textbf{Proof.} 
 Recall the following 2-cocycle condition. (See a paragraph after
 Definition \ref{df:2-cocycle}.)
 $$(\alpha_\pi\otimes \id)(U_{\rho,\sigma}^{\eta,a})
 U_{\pi,\eta}^{\xi,b} 
 =\sum_{\zeta,c,d}
 (U_{\pi,\rho}^{\zeta,c}\otimes 1_\sigma) U_{\zeta, \sigma}^{\xi,d} 
 V_{(\zeta,c,d),(\eta_{a,b})}.$$
 We put $\bar{\rho}=\sigma=\pi$, $\eta=\mathbf{1}$ (hence $\xi=\pi$), and multiply
 $U_{\pi,\bar{\pi}}^{\mathbf{1}}\otimes 1_\pi$ from the left on both sides. 
 Then we get
 the following.
  $$(U_{\pi,\bar{\pi}}^\mathbf{1*}\otimes
  1_\pi)\alpha_\pi\otimes \id\left(U_{\bar{\pi},\pi}^\mathbf{1}\right)=
  V_{\mathbf{1},\mathbf{1}}=
 (T_{\pi,\bar{\pi}}^\mathbf{1*}\otimes
  1_\pi)(1_\pi\otimes T_{\bar{\pi},\pi}^\mathbf{1})=\frac{1}{d\pi}.
 $$
 Since $U_{\pi_i,\bar{\pi}_j}^{\mathbf{1}}=\sum_k
 1/\sqrt{d\pi}U_{\pi_{ik},\bar{\pi}_{jk}}=1/\sqrt{d\pi}\tilde{U}_{\pi_i,\bar{\pi}_j}$, 
 we get the conclusion.
  \hfill$\Box$ \\

 \noindent
 \textbf{Proof of Lemma \ref{lem:reltwist}.} It is easy to see $\alpha$ is an action of $M$ on
 $H\otimes B(H_\pi)$.
 We verify that $\lambda_\pi$ implements $\alpha_\pi$. Then 
 \begin{eqnarray*}
  \left(\lambda_{\pi_{ij}}\alpha(a)v\right)(\rho)&=&
 \sum_{\sigma,e}U_{\rho,\pi_i}^{\sigma,e}\left(\alpha(a)v\right)(\sigma)
 T_{\rho,\pi_j}^{\sigma,e*} \\
 &=&
 \sum_{\sigma,e}U_{\rho,\pi_i}^{\sigma,e}\alpha_\sigma(a)v(\sigma)
 T_{\rho,\pi_j}^{\sigma,e*} \\
 &=&
 \sum_{\sigma,e,k}\alpha_\rho(\alpha_{\pi}(a)_{ik})
 U_{\rho,\pi_k}^{\sigma,e}v(\sigma)
 T_{\rho,\pi_j}^{\sigma,e*} \\
 &=&\sum_k\left(\alpha(\alpha_{\pi}(a)_{ik})\lambda_{\pi_{kj}}v\right)(\rho).
 \end{eqnarray*}
 holds. Therefore we have $\lambda_\pi(a\otimes 1_\pi)=\alpha_\pi(a)\lambda_\pi$ by
 identifying $\alpha(a)$ and $a$.

 We next compute $\lambda_\pi^{12}\lambda_\sigma^{13}$ as follows.
 \begin{eqnarray*}
  \left(\lambda_{\pi_{ij}}\lambda_{\rho_{kl}}v\right)(\xi)&=&
 \sum_{\sigma,a,\eta,b}
 U_{\xi,\pi_{i}}^{\sigma,a}U_{\sigma,\rho_k}^{\eta,b}v(\eta)T_{\sigma,\rho_l}^{\eta,b*}
 T_{\xi,\pi_j}^{\sigma,a*} \\
 &=&
 \sum_{\sigma,a,m,n,\eta,b}
 \alpha_\xi(U_{\pi_i,\rho_k}^{\sigma_m,a})U_{\xi,\sigma_m}^{\eta,b}
 v(\eta)T_{\xi,\sigma_n}^{\eta,b*} T_{\pi_j,\rho_l}^{\sigma_n,a*} \,\,\,
 (\mbox{by $2$-cocycle condition})\\
 &=&\sum_{\sigma,m,n,a}\left(\alpha(U_{\pi_i,\rho_k}^{\sigma_m,a})
 \lambda_{\sigma_{m,n}}T_{\pi_j,\rho_l}^{\sigma_n,a}v\right)(\xi).
 \end{eqnarray*}
 Hence we have
 $\lambda_{\pi}^{12}\lambda_\rho^{13}T_{\pi,\rho}^{\sigma,a}=
 U_{\pi,\rho}^{\sigma,a}\lambda_\sigma$.

 Finally, we verify that $\lambda_\pi$ is a unitary.
 One can easily to see
 $\sum_{k}\lambda_{\pi_{ki}}^*\lambda_{\pi_{kj}}=\delta_{i,j}$ (hence
 $\lambda_\pi^*\lambda_\pi=1$)
 from the definition of $\lambda_{\pi_{ij}}$ and
 $U_{\pi,\rho}^{\sigma,a*}U_{\pi,\rho}^{\xi,b}=\delta_{\sigma,\xi}\delta_{a,b}$. 
 Hence we only have to see
 $\lambda_\pi\lambda_\pi^*=1$. 

 To this end, we first show $\lambda_{\bar{\pi}_{ij}}^*=\sum_k 
  \tilde{U}_{\pi_{k},\bar{\pi}_i}^*\lambda_{\pi_{kj}}.$ 
 Since 
 $$\sum_{k}\lambda_{\pi_{ik}}\lambda_{\bar{\pi}_{jk}}=\sum_{\rho,l,m,a}
 U_{\pi_i,\bar{\pi}_j}^{\rho_l,a}\lambda_{\rho_{lm}}T_{\pi_k,\bar{\pi}_k}^{\rho_m,a}=
 \tilde{U}_{\pi_i,\bar{\pi}_j},$$ 
 we  have
 $\lambda_\pi{}^t\lambda_{\bar{\pi}}=\left(\tilde{U}_{\pi_i,\bar{\pi}_j}\right)_{i,j}.$ 
 Then we get
 $^t\lambda_{\bar{\pi}}=\lambda_\pi^*\left(\tilde{U}_{\pi_i,\bar{\pi}_j}\right)_{i,j}$. 
  Comparing matrix
  elements of both sides, we get $\lambda_{\bar{\pi}_{ij}}^*=\sum_k 
  \tilde{U}_{\pi_{k},\bar{\pi}_i}^*\lambda_{\pi_{kj}}.$

  Then we get 
  \begin{eqnarray*}
   \sum_k\lambda_{\pi_{ik}}\lambda_{\pi_{jk}}^*&=&\sum_{k,l}
  \lambda_{\pi_{ik}}\tilde{U}_{\bar{\pi}_l,\pi_j}^*\lambda_{\bar{\pi}_{lk}} \\
  &=& 
  \sum_{k,l,m}
  \alpha_\pi\left(\tilde{U}_{\bar{\pi}_l,\pi_j}^*\right)_{im}\lambda_{\pi_{mk}}
  \lambda_{\bar{\pi}_{lk}} \\
  &=& 
  \sum_{k,l,m,\xi,a,b}
  \alpha_\pi\left(\tilde{U}_{\bar{\pi}_l,\pi_j}^*\right)_{im}U_{\pi_m,\bar{\pi}_l}^{\xi_a}
  \lambda_{\xi_{ab}}T_{\pi_k,\bar{\pi}_k}^{\xi_b} \\
  &=&\sum_{k,l,m}
  \alpha_\pi\left(\tilde{U}_{\bar{\pi}_l,\pi_j}^*\right)_{im}\tilde{U}_{\pi_m,\bar{\pi}_l}\\
  &=&\delta_{i,j}
  \end{eqnarray*}
  by Lemma \ref{lem:u}, and $\lambda_\pi$ is indeed a unitary.

  \hfill$\Box$

  We construct a conditional expectation $E$ 
  from
   $M\rtimes_{\alpha,U}\hG$ onto $M$. Let $P$ be a projection from
   $H\otimes \ell^2(\hG) $ to $H\otimes B(H_\mathbf{1})\cong H$, and set 
  $E(x):=PxP^*$. Then $E$ is indeed a conditional expectation 
  from
   $M\rtimes_{\alpha,U}\hG$ onto $M$ with $E(\lambda_{\pi_{ij}})=\delta_{\mathbf{1},\pi}$.
  Then the following lemma can be easily verified as in the usual crossed product.

  \begin{lem}
  Every $a\in M\rtimes_{\alpha, U}\hG$ is 
  expressed uniquely as $a=\sum_{\pi,i,j}a_{\pi,i,j}\lambda_{\pi_{ij}}$,
   $a_{\pi,i,j}\in M$.  
  \end{lem}
  Here we only remark that a coefficient $a_{\pi,i,j}$ is given by
  $a_{\pi,i,j}=d\pi E(a\lambda_{\pi_{ij}}^*)$.

\ifx\undefined\bysame
\newcommand{\bysame}{\leavevmode\hbox to3em{\hrulefill}\,}
\fi

\end{document}